\documentclass[12pt]{amsart}
\usepackage{amsbsy}
\usepackage{graphicx,epsfig,subfigure}
\begin{document}

%
%
%

\newtheorem{theorem}{Theorem}
\newtheorem{proposition}[theorem]{Proposition}
\newtheorem{lemma}[theorem]{Lemma}
\newtheorem{corollary}[theorem]{Corollary}
\newtheorem{definition}[theorem]{Definition}
\newtheorem{remark}[theorem]{Remark}
\numberwithin{equation}{section} \numberwithin{theorem}{section}
\newcommand{\beq}{\begin{equation}}
\newcommand{\eeq}{\end{equation}}
\newcommand{\re}{{\mathbb R}}
\newcommand{\n}{\nabla}
\newcommand{\ren}{{\mathbb R}^N}
\newcommand{\iy}{\infty}
\newcommand{\pa}{\partial}
\newcommand{\ms}{\medskip\vskip-.1cm}
\newcommand{\mpb}{\medskip}
\newcommand{\ssk}{\smallskip}
\newcommand{\BB}{{\bf B}}
\newcommand{\Am}{{\bf A}_{2m}}
\newcommand{\bL}{\BB^*}
\newcommand{\bLs}{\BB}
\renewcommand{\a}{\alpha}
\renewcommand{\b}{\beta}
\newcommand{\g}{\gamma}
\newcommand{\ka}{\kappa}
\newcommand{\G}{\Gamma}
\renewcommand{\d}{\delta}
\newcommand{\D}{\Delta}
\newcommand{\e}{\varepsilon}
\newcommand{\vp}{\varphi}
\renewcommand{\l}{\lambda}
\renewcommand{\o}{\omega}
\renewcommand{\O}{\Omega}
\newcommand{\s}{\sigma}
\renewcommand{\t}{\tau}
\renewcommand{\th}{\theta}
\newcommand{\z}{\zeta}
\newcommand{\wx}{\widetilde x}
\newcommand{\wt}{\widetilde t}
\newcommand{\noi}{\noindent}
\newcommand{\lb}{  (}
\newcommand{\rb}{  )}
\newcommand{\lsb}{  [}
\newcommand{\rsb}{  ]}
\newcommand{\lab}{  \langle}
\newcommand{\rab}{  \rangle }
\newcommand{\gap}{\vskip .5cm}
\newcommand{\bz}{\bar{z}}
\newcommand{\bg}{\bar{g}}
\newcommand{\Ba}{\bar{a}}
\newcommand{\bt}{\bar{\th}}
\def\com#1{\fbox{\parbox{6in}{\texttt{#1}}}}

\title
{\bf Self-similar blow-up in parabolic equations of
Monge--Amp\`ere type}

\author {C.J.~Budd and V.A.~Galaktionov}

\address{Department of Mathematical Sciences, University of Bath,
 Bath BA2 7AY, UK}
\email{cjb@maths.bath.ac.uk}

\address{Department of Mathematical Sciences, University of Bath,
 Bath BA2 7AY, UK}
\email{vag@maths.bath.ac.uk}

\date{\today}



\begin{abstract}

We use techniques from reaction-diffusion theory
to study the blow-up and existence of solutions of the parabolic Monge--Amp\`ere equation
with power source,  with the following basic 2D model
 \beq
 \label{0.1}
 u_t = -|D^2 u| + |u|^{p-1}u \quad \mbox{in} \quad \re^2 \times \re_+,
   \eeq
   where in two-dimensions
$| D^2 u|= u_{xx} u_{yy}- (u_{xy})^2$ and  $p > 1$ is a fixed
   exponent.
   For a class of ``dominated concave"
   and compactly supported radial initial data
   $u_0(x) \ge 0$, the Cauchy problem is shown to
   be locally well-posed and to exhibit finite time blow-up that
   is described by similarity solutions. For $p \in (1, 2]$,
   similarity solutions, containing domains of concavity and convexity,
   are shown to be
    compactly supported and
   correspond to surfaces with flat sides that persist until the
   blow-up time. The case $p > 2$ leads
to single-point blow-up. Numerical computations of blow-up solutions
without radial symmetry are also presented.

   The  parabolic analogy of  (\ref{0.1}) in 3D for which $|D^2 u|$ is
   a  cubic operator is
   $$
    u_t = |D^2 u| + |u|^{p-1}u \quad \mbox{in} \quad \re^3 \times \re_+,
   $$
   and is shown to
    admit a wider
   set of (oscillatory) self-similar blow-up patterns.
Regional  self-similar blow-up in a cubic radial model related to
the fourth-order M-A equation of the type
 $$
 u_t = -|D^4 u| + u^3
  \quad \mbox{in} \quad \re^2 \times \re_+,
   $$
where the cubic operator $|D^4 u|$ is the  catalecticant $3 \times
3$ determinant, is also briefly discussed.

This is an earlier extended version of \cite{BGMAIMA}, where, in
particular, we present a survey on various M-A models; see
Appendix A.

\end{abstract}

\maketitle

\section{Introduction:
  our basic parabolic Monge--Amp\`ere equations with blow-up}



\subsection{Outline}

\noindent Fully nonlinear parabolic partial differential equations
with spatial Monge--Amp\`ere (M-A) operators arise in many
problems related to optimal transport and geometric flows
\cite{Cull}, image registration \cite{HakT}, adaptive mesh generation
\cite{Delz, BuHu}, the evolution of vorticity in meteorological systems
\cite{Cull}, the semi-geostrophic equations of meteorology, as
well as being extensively studied in the analysis literature; see
Taylor   \cite[Ch.~14,15]{Tay}, Gilbarg--Trudinger
\cite[Ch.~17]{GilTr}, Guti\'errez \cite{Gut01}, and
Trudinger--Wang \cite{TrudW08}, as a most recent reference. To
describe such equations, we consider a given function $u \in
C^2(\ren)$, for which $D^2u$ denotes the corresponding $N \times
N$ {\em Hessian} matrix  $D^2 u=\|u_{x_ix_j}\|$, so that in
two-dimensions, $d=2$,
 \beq
  | D^2 u| \equiv {\rm det}\, D^2 u = u_{xx} u_{yy}-(u_{xy})^2.
\eeq
  Similarly, in three dimensions
 \beq
 \label{M133}
 |D^2 u|= \left|
  \begin{matrix}
  u_{xx} \,\, u_{xy} \,\, u_{xz}\\
  u_{xy} \,\, u_{yy} \,\, u_{yz}\\
  u_{xz} \,\, u_{yz} \,\, u_{zz}
  \end{matrix}
  \right|.
  \eeq
A general  parabolic {\em Monge--Amp\`ere} (M-A) {\em equation}
with a nonlinear source term, then takes the  form
 \beq
  \label{MA.1}
   \mbox{$
  u_t = g\bigl({\rm det} D^2  u\bigr) + h(x,u,D_x u) \quad \mbox{in} \,\,\, \ren \times
  \re_+
   $}
 \eeq
with proper initial data $u(0,x) = u_0(x)$. Such PDEs with various
nonlinear operators $g(\cdot)$ and $h(\cdot)$  have a number of
important applications.

\vspace{0.1in}

\noindent The~origin of such fully nonlinear M-A equations dates
back to Monge's paper \cite{Monge81} in 1781,\index{Monge, G.} in
which Monge proposed a civil-engineering problem of moving a mass
of earth from one configuration to another in the most economical
way. This problem has been further studied by
 Appel \cite{Appel87} and L.V. Kantorovich \cite{Kant42, Kant48};
  see
references and a survey in \cite{Feyel03}. Other key problems and
M-A applications include: logarithmic Gauss and Hessian  curvature
flows, the Minkowski problem (1897) and the Weyl problem (with
Calabi's related conjecture in complex geometry), etc.

\vspace{0.1in}

 For increasing functions $g(s)$, the equation (\ref{MA.1}) is
parabolic if $D^2u(\cdot,t)$ remains positive definite for $t>0$,
assuming that $D^2u_0>0$ and local-in-time solutions exist
\cite[p.~320]{Lad68}. Provided that $g(s)$ does not grow too
rapidly, for example if $$ \mbox{$
 g(s) = \ln s, \,\,\, g(s) = - \frac 1s, \,\,\,\mbox{and} \,\,\,
g(s) = s^{\frac 1N} \quad \mbox{for} \,\,\, s> 0, $}$$ it is known
\cite{Iv94, Gut01, Chou00} that the solutions of M-A exist for all
time.

\noindent In general, however, the  questions of local solvability
and regularity for M-A equations even in 2D such as the system
 \beq
 \label{Kh1}
 (u_{xx}+a)(u_{yy}+c)-(u_{xy}+b)^2=f
 \eeq
in the hyperbolic ($f<0$) and mixed type ($f$ of changing sign)
 are difficult, and there are some counterexamples concerning
 these basic theoretical problems and concepts;
 see \cite{Khuri07}. Note that classification problem for the M-A
 equations such as (\ref{Kh1}) on finding their simplest form was
 already posed and partially solved by Sophus Lie in 1872-74
 \cite{Lie72}; see details and recent results in \cite{Kush08}.

\vspace{0.1in}

\noindent For other  functions $g(s)$ in (\ref{MA.1}) with a
faster growth as $s \to \infty$ and for certain nonlinear source
terms $h$, solutions which are locally regular may evolve to
blow-up in a finite time $T$. This  gives  special singular
asymptotic patterns,
  which can
  also be of interest in some geometric applications; on singular patterns for M-A flows, see \cite{DH99, DL03}
The analysis literature currently suffers from a lack of
understanding about the formation of such singularities in the
fully nonlinear M-A equation despite their relevance to such
problems as front  formation in meteorology \cite{Cull}. This
paper aims to make a start at studying such blow-up singular
behaviour by using techniques derived from studying
reaction-diffusion equations to look at some special parabolic M-A
problems, which lead to the finite time formation of
singularities. In particular,
 we consider (\ref{M1}), (\ref{M133}), and some other higher-order PDEs
  as formal basic
 equations demonstrating that M-A models can exhibit several common
features of blow-up, which have been previously observed in PDEs
with classic reaction-diffusion, porous medium,  and the
$p$-Laplacian operators. Indeed, we will approach the study of M-A
type operators by developing the related theory of the
$p$-Laplacian operator. Our interest in this paper will be an
understanding of the various forms of singularity that can arise
in some  parabolic  Monge--Amp\`ere models with a {\em polynomial
source term} as well as extending the general existence theory for
such problems.

Further discussion of a variety of M-A parabolic, elliptic, and
hyperbolic models together with basic regularity and singularity
results continues  in Appendix A after the list of references.

\subsection{Model equations and results}

\noindent {\bf Model 1} Our  first model fully nonlinear PDE is given by
 \beq
 \label{M1}
  u_t = (-1)^{d-1}| D^2 u| + |u|^{p-1}u \quad \mbox{in}
\quad \re^d \times \re_+,
   \eeq
    where $p>1$ is a given constant. Such models are natural counterparts of the porous
    medium equation  with
 reaction/absorption, and  of thin film (or Cahn--Hilliard-type, $n=0$)
 models,
  \beq
  \label{RD.991}
  u_t = \Delta  u^m \pm u^p \quad
  \mbox{and} \quad  u_t= - \n \cdot (|u|^n \n \Delta  u) \pm \Delta
  |u|^{p-1}u.
  \eeq
Our principle interest lies in the study of those solutions which have
large isolated spatial maxima tending towards singularities forming in the finite time $T$.
Such solutions are locally concave close to the peak.
The choice of sign of the principal operator in (\ref{M1})
ensures local well-posedness (local parabolicity) of the partial
differential equation in such neighborhoods. Significantly, the
existence theory for such locally concave solutions is rather different
from the usual theory of the  M-A operator, which is restricted
to globally convex solutions, and we will look at it detail
in Section 4. The
initial data $u_0(x) \ge 0$ is
 assumed to be bell-shaped (this preserves ``dominated concavity")
and sometimes compactly
  supported.
 We firstly study radially symmetric solutions in two and three spatial dimensions,
 and will show analytically, by extending the theory of
$p$-Laplacian operators, and demonstrate numerically, that whilst
the Cauchy problem
 is locally well-posed, and admits a unique radially symmetric
weak solution, certain of these solutions become
singular with finite-time blow-up. We will also find a set
of self-similar blow-up patterns,
corresponding to single-point blow-up if $p > d$,
 regional blow-up for $p=d$, and
to global blow-up for $p \in(1,d)$. The stability of these will
be investigated numerically, and we will find that monotone self-similar
blow-up profiles appear to be globally stable.

\vspace{0.1in}

\noindent We will also present some analytic and numerical
results for the time evolution of non radially symmetric
solutions in two dimensions.
These computations will give some evidence to conclude that stable
 non-radially symmetric blow-up solution profiles also exist,
 though this leads to a number of difficult open mathematical problems.

\vspace{0.1in}

\noindent {\bf Model 2} As a second model equation, we will look at higher-order
fully nonlinear  M-A spatial operators associated with the
equations of the form
\begin{equation}
{ u_t = (-1)^{d-1} |D^4u| + |u|^{d}u \quad \mbox{in} \quad \re^d
\times \re_+, }
 \label{1.cjb}
\end{equation}
where $|D^4 u|$ is the determinant of the 4th derivative matrix of
$u$. We will find similar results on the blow-up profiles to those
for the second-order operator. A principal feature of compactly
supported solutions to (\ref{1.cjb}) is that these are infinitely
oscillatory and changing sign at the interfaces, and this property
persists until blow-up time.

\vspace{0.2in}

\noindent The layout of the remainder of this paper is as follows.
In Section 2 we look at single-point blow-up, regional blow-up,
and the global one of the radially symmetric solutions of the
polynomial M-A equation in two and three spatial dimensions.  We
will combine both an analytical and a numerical study to determine
the form and stability of the self-similar blow-up solutions. In
Section 3, we will extend this analysis to look at solutions,
which do not have radial symmetry, and will give numerical
evidence for the existence of stable blow-up profiles in non
radial geometries. In Section 4, we will look at more general
properties associated with M-A type flows, in particular the
existence of various conservation laws. Finally, in Section 5 we
will study the forms of  the blow-up behaviour for the equations
with higher-order operators as in (\ref{1.cjb}).

\section{Parabolic M-A equations in $\re^2$: blow-up in radial geometry}
 \label{S2}

\subsection{Radially symmetric solutions: first results on blow-up}

\noindent The Hessian operator $|D^2 u|$ restricted to radially
symmetric solutions in $\re^{d}$ takes the form of a
non-autonomous version of the $p$-Laplacian operator. Namely,
 for solutions $u=u(r,t)$,
  with the single spatial variable $r=|(x,y,...)| >0$, equation (\ref{M1})
takes the form
 \beq
 \label{2.1}
  \mbox{$
 u_t = (-1)^{d-1} \frac{1}{r^{d-1}} \, (u_r)^{d-1} u_{rr} + |u|^{p-1}u \quad \mbox{in}
 \quad \re_+ \times \re_+,
  $}
  \eeq
and then for $r=0$ we have the symmetry condition $$u_r = 0.$$

\noindent In this section, we shall consider the nature of the
blow-up solutions for this problem and will identify different
classes (single-point, regional and global) of self-similar radial
solutions, giving some numerical evidence for their stability.
However, we note at this stage (and will establish in the next
section) that (possibly stable) non-radially symmetric blow-up
solutions of the underlying PDE also exist. One can see that
(\ref{2.1}) implies that the equation is (at least, degenerate)
parabolic if
 \beq
 \label{par1}
 (-1)^{d-1} (u_r)^{d-1} \ge 0 \quad
  \Longrightarrow \quad
 \left\{
  \begin{matrix} u_r \le 0 \quad \mbox{for even $d$},\qquad\,\,\, \ssk \\
 u_r\,\,\, \mbox{is arbitrary for odd $d$}.
 \end{matrix}
  \right.
 \eeq
For the local well-posedness of the above M-A flow, (\ref{par1})
is always assumed.
 In all the cases, the differential operator $(-1)^{d-1} \frac{1}{r^{d-1}}
\, (u_r)^{d-1} u_{rr}$ is regular in the class of monotone
decreasing, sufficiently smooth, and strictly concave at the
origin functions, so that the local well-posedness of (\ref{2.1})
is guaranteed for the initial data satisfies the regularity and
monotonicity  constraints
 \beq
 \label{2.2}
 u(r,0)=u_0(r) \ge 0, \quad u_0 \in C^1([0,\infty)), \quad u_0'(0)=0, \quad
 u_0'(r) \le 0 \,\,\, \mbox{for} \,\,\, r>0.
  \eeq
  The corresponding $p$-Laplacian counterpart of (\ref{2.1}) is
 then
 \beq
 \label{2.3}
 \mbox{$
 u_t = \frac{1}{r^{d-1}} \, |u_r|^{d-1} u_{rr} + |u|^{p-1}u \quad \mbox{in}
 \quad \re_+ \times \re_+,
  $}
  \eeq
  which is locally well-posed by parabolic regularity theory; see
  e.g., \cite{DB}.
   By the Maximum Principle (MP), the assumptions
  in (\ref{2.2}) guarantee that the solution $u=u(r,t)$ satisfies
  the monotonicity condition
   \beq
   \label{2.4}
    u_r(r,t) \le 0.
    \eeq
\noindent Therefore, equations (\ref{2.1}) and (\ref{2.3})
coincide in this class of monotone solutions. Note that for local
well-posedness, we do not need any {\em concavity-type}
assumptions that are usual for standard M-A flows.




\vspace{0.1in}

\noindent The phenomenon of blow-up for the solutions of the model (\ref{2.3})
can be studied by using techniques derived from the
study of reaction diffusion equations
(see
\cite[Ch.~4]{SGKM}).  By a comparison of the solution with sub-and super-solutions of the same equation of
self-similar form, we can show
  that, for nonnegative solutions $u$,
   there exists a
  critical Fujita exponent
  \beq
  \label{Fuj1}
  p_0=d+2
    \quad \mbox{such that:}
  \eeq

 (i) for $p \in (1,p_0]$, any $u(x,t) \not \equiv 0$ blows up in
 finite time, and

 (ii) in the supercritical range $p>p_0=d+2$, solutions blow-up
 for large enough data, while for small ones, the solutions are
 global in time.

 \noindent The proof of blow-up in the critical case $p=p_0$ is most
 delicate and demands a monotonicity/asymptotic rescaled construction; see
 e.g., \cite{GalFuj, GPohTFE}.

\noindent For the remainder of this paper we shall always assume that the initial data
are such that blow-up always occurs.

\subsection{Blow-up similarity solutions}

\vspace{0.1in}

\noindent The M-A equation with the $|u|^{p-1}u$ source term is
invariant under the scaling group
  $$
  t \to \lambda t, \quad r \to
\lambda^{\frac{p-d}{2d(p-1)}}r, \quad u \to \lambda^{- \frac
1{p-1}}u.
  $$
 Accordingly, a self-similar blow-up profile (with
blow-up at the origin $r=0$, which is assumed to belong to the
blow-up set) is described by the following solutions:
 \beq
 \label{M5}
 \mbox{$
 u_S(r,t)= (T-t)^{-\frac 1{p-1}} f(z), \quad z= \frac{r}{(T-t)^\b}, \quad
 \b= \frac{p-d}{2d(p-1)}.
  $}
  \eeq
  Here $f \ge 0$ is a solution of the following ordinary
  differential equation,
  \beq
  \label{M6a}
   \mbox{$
  \frac{1}{z^{d-1}} \, (-1)^{d-1}(f')^{d-1} f'' - \b f'z- \frac {1}{p-1}\, f + |f|^{p-1}f=0,
   \,\,\, f'(z)
    \le
  0,
 \,\, f'(0)=f(+\infty)=0.
  $}
  \eeq
\noindent The condition on $f(+\infty)$ (and the consequent
requirement that the solutions of (\ref{M6a}) should decay to zero
as $z \to \infty$) is necessary to ensure that the self-similar
solutions correspond to solutions of the original Cauchy problem
of the PDE. \noindent In the case of monotone decreasing solutions
with $f'(z) \le 0$ the equation (\ref{M6a}) becomes \beq
  \label{M6}
   \mbox{$
  \frac{1}{z^{d-1}} \, |f'|^{d-1} f'' - \b f'z- \frac {1}{p-1}\, f + |f|^{p-1}f=0,  \,\,\, f'(z) \le
  0\,\,\,
 \mbox{in} \,\,\, \re_+; \quad f'(0)=f(+\infty)=0.
  $}
  \eeq

\noindent In particular, in the case of $d=2$ to be studied in
greater detail, we have to require the monotonicity assumption to
allow the construction of smooth solutions. In the case of $d=3$,
we can relax this assumption, leading to a richer class of
(possibly oscillatory) solutions. When considered as an initial
value problem, the solutions of (\ref{M6a}) lose regularity at the
degeneracy $f'=0$ (except the origin $r=0$) and are not twice
differentiable at such points. However the existence of weak
solutions with reduced regularity is guaranteed by the standard
theory of the p-Laplacian operator, and these questions are
standard in parabolic theory; see \cite{DB} and
\cite[Ch.~2]{AMGV}.  Note that the  ODE (\ref{M6}) has two
constant equilibria given by
  \beq
  \label{eq11}
 \pm f_0 = \pm (p-1)^{-\frac 1{p-1}}
   \eeq
and that solutions close to these equilibria can be oscillatory,
which is a crucial property to be properly treated and used.

\vspace{0.1in}

\noindent If $p > 1$, then each solution of the form (\ref{M5})
blows up in finite time, however the nature of the scaling is very
different in the three cases of $1 < p < d$, $p = d$ and $p > d$,
corresponding to global (blow-up over the whole of $\re^d$),
regional (blow-up over a sub-set of $\re^d$ with non-zero
measure),
 and single-point
blow-up respectively (zero measure blow-up set in general).
Indeed, we will show that for $p \in(1,d]$ the solutions $f(z)$ of
the ODE are compactly supported, while for $p>d$ they are strictly
positive.


\subsection{Regional blow-up when $p=d$}


\noindent We begin with the case $p=d$, where, according to
(\ref{M5}), the ODE becomes autonomous and $z = r$, so that
 \beq
 \label{M7}
  \mbox{$
  \frac{1}{r^{d-1}} \, |f'|^{d-1} f'' - f  + |f|^{p-1}f=0,  \,\,\, f'(r) \le
  0\,\,\,
 \mbox{in} \,\,\, \re_+; \quad f'(0)=f(+\infty)=0.
  $}
 \eeq
This problem falls into the scope of the well-known blow-up
analysis for quasilinear reaction-diffusion equations,
\cite[Ch.~4]{SGKM}.

\begin{proposition}
 \label{Pr.1}
 The problem $(\ref{M7})$ has a non-trivial, monotone, compactly supported solution
 $F_0(z) \ge 0$ such that $F_0(0)>1$. The support of
$F_{0}(z)$ is given by $[0,L_{S}]$ with the asymptotic behaviour
near the interface:  as $z \to L_S$,
  \beq
  \label{as1N}
   \mbox{$
F_{0}(z) = A (L_{S} - z)_+^{\frac{d+1}{d-1}}(1+o(1)), \quad
\mbox{where} \quad A = \big[
\frac{(d-1)^{d+1}}{2(d+1)^d}\big]^{\frac 1{d-1}}.
 $}
 \eeq
  \end{proposition}

  \noi{\em Proof.} The result follows the lines of the ODE
  analysis in \cite[pp.~183-189]{SGKM} and uses a shooting approach
in which (\ref{M6a}) is considered as an IVP with shooting
parameter $f(0)$. If $f(0)$ is too large then the solutions of the
IVP diverge to $-\infty$ and if it is not large enough then the
solution has an oscillation about $f_0$ in a manner to be
described in more detail below. The self-similar solution occurs
at the point of transition between these two forms of behaviour.
$\qed$

\vspace{0.1in}

\begin{figure}
\centering
\includegraphics[scale=0.6]{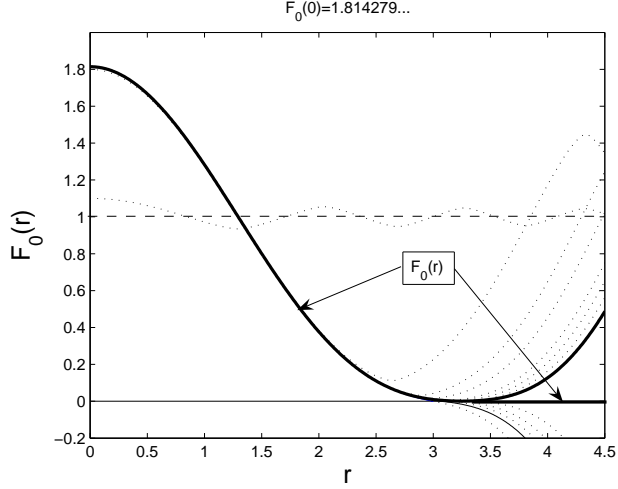}
\caption{\small $p=d=2$: the similarity profile $F_0$ obtained by
shooting in the ODE (\ref{M7}) from the origin $r=0$ with the
parameter of shooting given by $f(0)>0$.}
 \label{F1}
\end{figure}

\noindent The form of the proof leads to a {\em numerical method based on shooting} for constructing
an approximation to
the function $F_{0}$. To do this we
specify the value of $F_{0}(0)$, take $F'_{0}(0) = 0$ and solve
(\ref{M7}) as an initial value problem in $r$ using an accurate
numerical method (typically a variable step BDF method). In this
numerical calculation the term $|f'|$ is replaced by $\sqrt{\e^2 +
(f')^2}$ with $\e = 10^{-5}$. This allows the numerical method to
cope with the loss of regularity when $f'=0$. The value of $f(0)$
is then steadily increased from $f_0$ and the transition between
oscillatory and divergent behaviour determined by bisection. At
this particular value  $f(0) \equiv F_0(0)$ the solution
approaches zero as $z \to L_{S}$ and for $z
>  L_{S}$ we have $F_{S}(0) \equiv 0$. Unlike studies of reaction-diffusion equations,
 the proof of the uniqueness of the solution
of (\ref{M7}) is not straightforward. However, the numerical
calculations strongly indicate the conjecture that such a
compactly supported monotone decreasing profile $F_0(z)$ is indeed
unique. In the case of $d=2$,
 its support is
 \beq
 \label{M8}
 L_S= 3.26...\, ,  \quad \mbox{with} \quad F_0(0)=1.814279... \, .
  \eeq
Similarly, in the case of $d=3$, \beq
 \label{M8cjb}
 L_S= 2.303...\, ,  \quad \mbox{with} \quad F_0(0)=1.366... \, .
  \eeq
 In Figure \ref{F1}, we show the resulting profile (bold) in the case of $p=d=2$,
obtained by shooting as described above, together
with an oscillatory solution of the IVP when $f(0)  \approx f_0$ and some
nearby divergent solutions.

\begin{figure}
\centering
\includegraphics[scale=0.6]{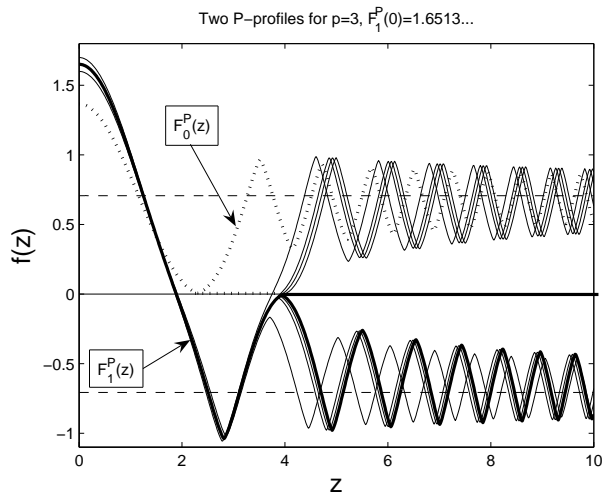}
\caption{\small $p=d=3$: the similarity profiles $F_0$ (dotted)
and $F_1^P$ (bold) and some nearby oscillatory solutions, obtained
by shooting in the ODE (\ref{M7}) from the origin with the
parameter of shooting $f(0)>0$.}
 \label{FSh2}
\end{figure}

\noindent In the case of $d=3$, we can relax the monotonicity
requirement on the function $f'(z)$. In this case there exists a countable
set $\{F_k^P\}$ of compactly supported profiles that change sign
precisely $k$ times for any $k=0,1,2,... \, .$ A numerical
shooting calculation of both $F_0(z)$ (dotted) and of  $F_1^P(z)$
(bold), together with some nearby oscillatory solutions, is shown
in Figure \ref{FSh2}. This figure both explains how these further
profiles can be obtained numerically and indicates how their
existence can be justified rigorously along the lines of the
proofs in \cite{BuGa, GPos}. Here, we have
 $$
L_S^{(1)}= 3.95...\, ,  \quad \mbox{with} \quad F_1^P(0)=1.6513...
\, .
 $$

\vspace{0.1in}

\noindent It follows immediately that the variable separable
solution given by
 \beq
 \label{M9}
  \mbox{$
 u_{S}(r,t) = \frac 1{T-t}\, F_0(r)
  $}
  \eeq
 describes {\em regional blow-up}  which is localized in the disc/ball $\{r \le L_S\}$.

\vspace{0.1in}

\begin{figure}
\centering
\includegraphics[scale=0.4]{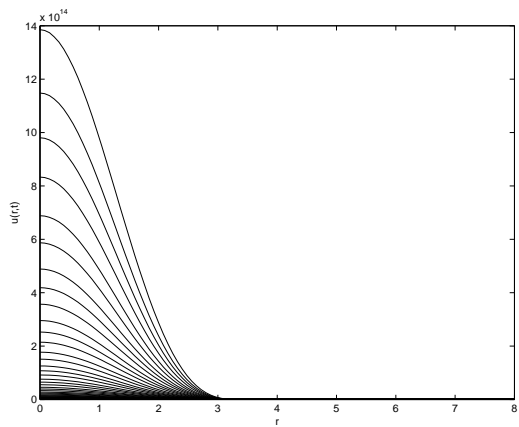}
\caption{\small $p=d=2$: Regional blow-up
of the function $u(r,t)$ with $\alpha = 0.1$.}
 \label{CJBF1}
\end{figure}

\begin{figure}
\centering
\includegraphics[scale=0.4]{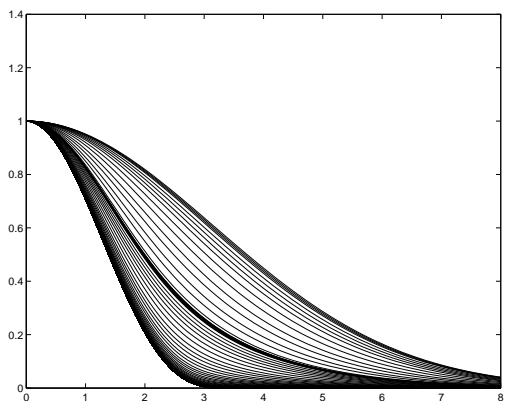}
\caption{\small $p=d=2$: Regional blow-up
of the scaled function $u(r,t)/u(0,t)$
with $\alpha = 0.1$ showing convergence to the
self-similar profile with support $[0,3.26]$ from above.}
 \label{CJBF2}
\end{figure}

\noindent We now make a further numerical calculation to investigate the stability of such a blow-up
profile (in the restricted class of radially symmetric solutions). To do this we use a semi-discrete
numerical method in which we discretise the Monge--Amp\`ere
spatial operator on a fine spatial mesh. This leads to a set of
ordinary differential equations for the discrete approximation to
the solution $u(r,t)$. These (stiff ordinary differential)
equations are then solved using an accurate variable order BDF
method. In this calculation we substitute a spatial domain $r \in
[0,L]$ for the infinite interval and impose a Neumann condition at
the boundaries $r = 0$ and $r = L$. For a calculation with $p = d
= 2$ we take $L = 8$ and use a spatial discretisation step size of
$L/1000$. We present in the following figures some calculation
showing the evolution of $u(r,t)$ and the scaled function
$u(r,t)/u(0,t)$ taking as initial data
   $$
   u(r,0) = 10\, {\mathrm
e}^{- \alpha x^{2}/2}.
  $$
   We consider two values of $\alpha$ to
give profiles which lie above and below the self-similar solution.
With $\alpha = 0.1$ the solutions initially lie above the
self-similar solution and in this case we see clear evidence in
Figures~\ref{CJBF1} and \ref{CJBF2} for evolution towards
self-similar regional blow-up with a blow-up time of $T \approx
0.1099$.

\begin{figure}
\centering
\includegraphics[scale=0.3]{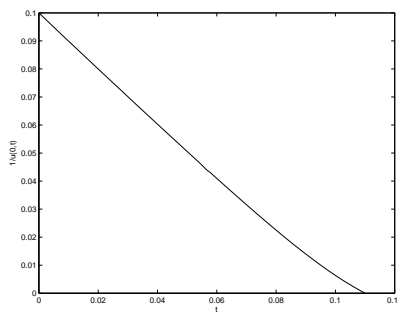}
\caption{\small $p=d=2$: Evolution of $1/u(0,t)$ with $\alpha =
0.1$ showing that $u(0,t)$ increases towards infinity. Note the
linear behaviour $\sim T-t$  close to the blow-up time.}
 \label{CJBF3}
\end{figure}

\noindent We also plot in  Figure~\ref{CJBF3} the value of $m(t) =
1/u(0,t)$. For small values of $m$ this figure is very  close to
linear, and a linear fit gives
 $$
  \mbox{$
 m(t) \approx -0.5553 t + 0.0610,
\quad \mbox{so that} \quad u(0,t) \approx \frac {1.805}{0.1099 -
t},
 $}
 $$
which is in good agreement with the earlier calculation of the
self-similar profile.

\noindent For comparison we now take initial data $\alpha = 10$.
In this case the value of $u(0,t)$ initially decreases, and then
increases with a blow-up time of $T \approx 1.4371$.
Asymptotically the blow-up is almost identical to that observed
earlier, however we can see clearly that this time the function
$u(r,t)/u(0,t)$ approaches the regional self-similar profile from
below.

\begin{figure}
\centering
\includegraphics[scale=0.5]{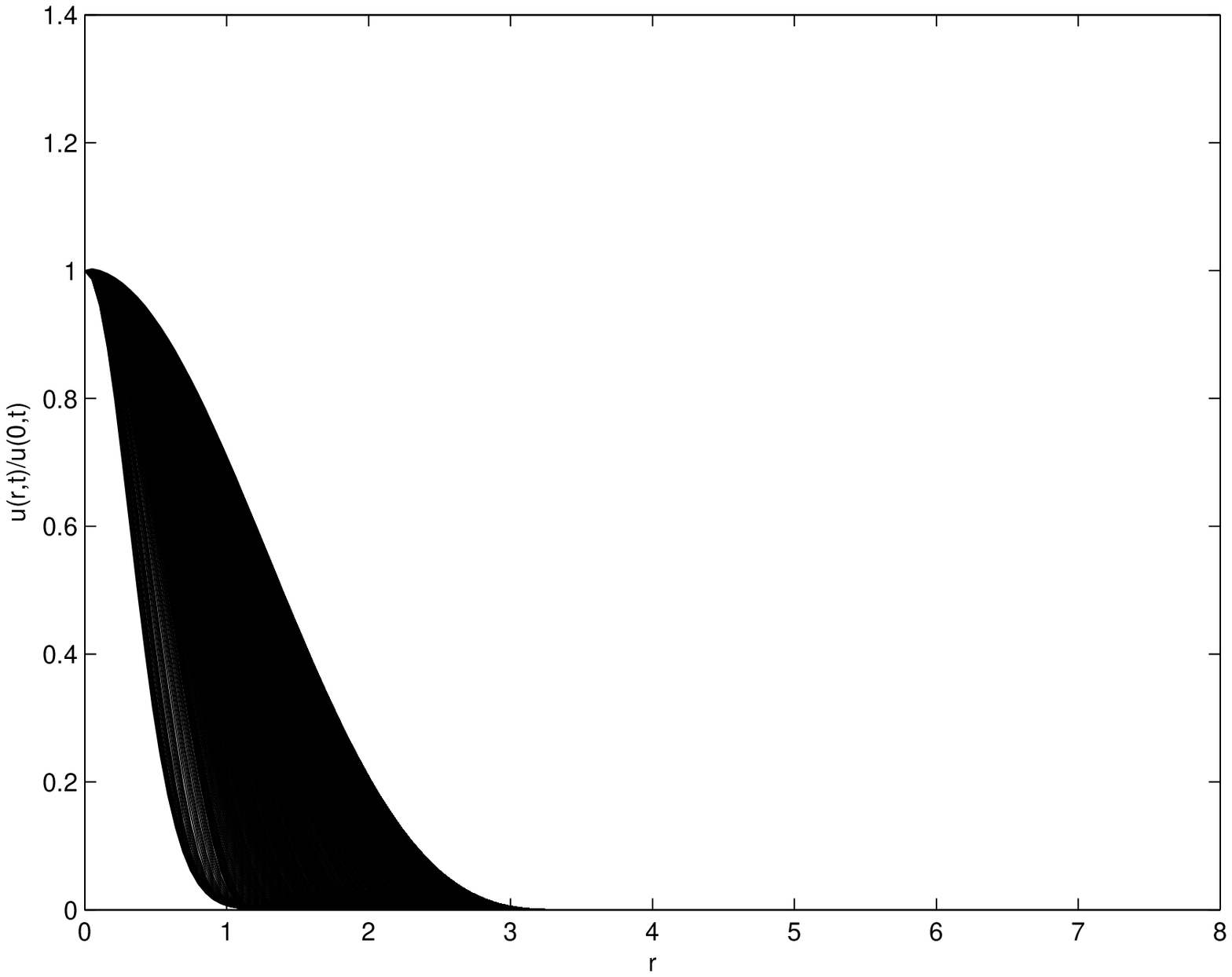}
\caption{\small $p=d=2$: Regional blow-up
of the scaled function $u(r,t)/u(0,t)$
with $\alpha = 10$ showing convergence to the
self-similar profile with support $[0,3.26]$ from below.}
 \label{CJBF5}
\end{figure}

\begin{figure}
\centering
\includegraphics[scale=0.3]{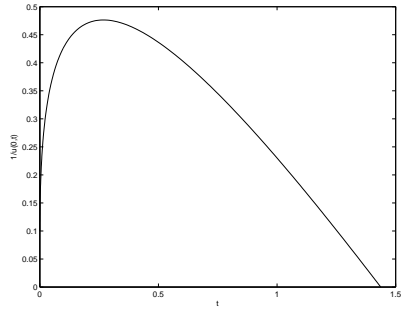}
\caption{\small $p=d=2$: Evolution of $1/u(0,t)$ with $\alpha =
10$ showing that $u(0,t)$ initially decreases before tending
towards infinity. Note the linear behaviour  $\sim T-t$ close to
the blow-up time.}
 \label{CJBF6}
\end{figure}

\noindent Very similar figures arise in the case of $d=3$
indicating that the regional self-similar solution
is also stable in this case.




\subsection{Single point blow-up for $p>d$: P and Q profiles}


\noindent If $p > d$ then $\beta > 0$ and single point blow-up occurs at the origin.
The self-similar blow-up profiles can then take various forms. Initially
we consider the monotone profiles for general $d$.

\begin{proposition}
 \label{PrLS}
 For any $p>d$,  the ordinary differential equation problem $(\ref{M6})$
 admits two
  strictly positive
 solutions $F_0^P(z) > 0$ and $F_0^Q(z) > 0$, which each satisfy the
asymptotic expansion
   \beq
   \label{C01}
 F_0(z) = C_0 z^{-\frac{2d}{p-d}}(1+o(1))
\,\,\, \mbox{as} \,\,\, z \to + \infty \quad (C_0>0),
   \eeq
and for which
 {\rm (i)} $F_0^P(z)$ is strictly monotone decreasing, $F_0'(z) < 0$
 for $z>0$,
with
  \beq
  \label{s2}
  F_0(0)> f_0, \quad \mbox{and}
\eeq

   {\rm (ii)} $F_0^Q(z)\equiv f_0$ on some interval $z \in [0,a_0]$
$($i.e., it has flat sides in this ball$)$,
 has the behaviour close to the interface given by
  \beq
  \label{ff11N}
 \mbox{$
F_0(z)= f_0- \frac 12 \, \b a_0^2(z-a_0)^2_+(1+o(1)) \quad
\mbox{as} \quad z \to a_0,
 $}
 \eeq
    and
is strictly monotone decreasing and smooth for $z > a_0$.

\end{proposition}

\begin{figure}
\centering
\includegraphics[scale=0.6]{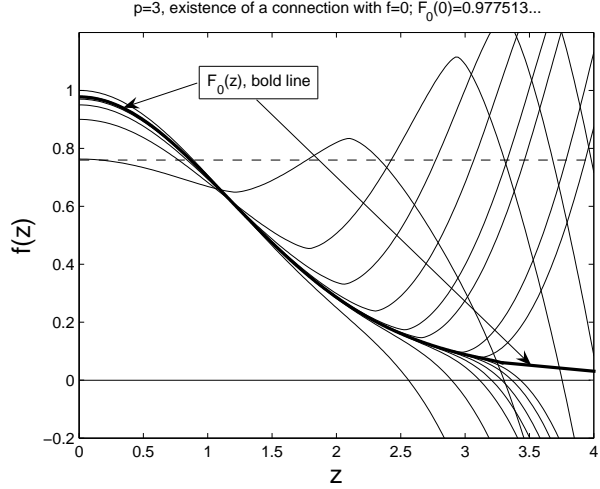}
\caption{\small Similarity P-profile $F_0$ for $d=2$ and $p=3$ in
(\ref{M6}) obtained by shooting from the origin $y=0$ with the
parameter of shooting $F_{0}(0)>0$.}
 \label{F2}
\end{figure}

\noindent The main ingredients of the proof of both P-type
profiles of the form (i)
and Q-type profiles of type (ii)
are explained in \cite{GPos, BuGa}. Both solutions are again
constructed by shooting, with $f(0)$ being the shooting parameter
for the P-type profiles and $a_0$ for the Q-type profiles. In
Figure \ref{F2}, we show the similarity profile of the P-type
solution $F_0(z)$ (bold) for the case of $d=2$ and $p=3$ together
with some nearby divergent solutions. In this case we find that
$F_{0}(0) = 0.9751...$\, . Numerically there is strong evidence
for uniqueness of this solution, but a proof of this uniqueness
result is open.

\vspace{0.1in}

\begin{figure}
\centering
\includegraphics[scale=0.6]{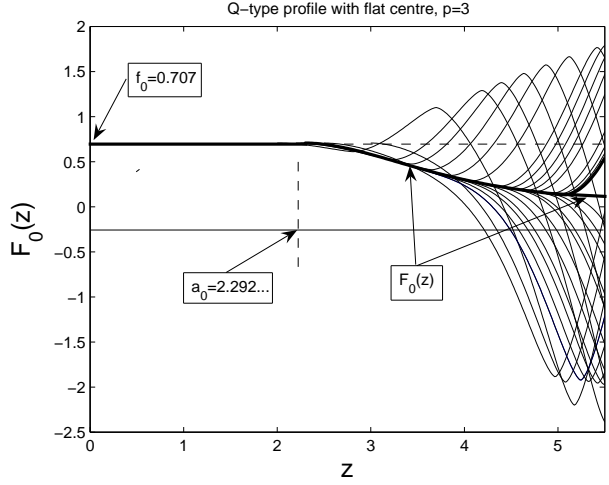}  
\caption{\small Similarity Q-profile $F_0$ for $d=2$, $p=3$ in
(\ref{M6}) obtained by shooting from $z=a_0>0$, which is the
parameter of shooting.}
 \label{F2Q}
\end{figure}

\noindent Similarly, in Figure \ref{F2Q}, we show
the results of shooting to find
 the Q-type profile when
for $d=2$ and $p=3$ and we obtain numerically that
$a_0=2.292...$\,.

\vspace{0.2in}

\noindent In the case of $d=3$, we may extend these results to
construct countable families of non-monotone P and Q-type
solutions described as follows:

\begin{proposition}
 \label{PrLS3}
 If $d=3$ then for any $p>3$, problem $(\ref{M6a})$
 admits the following two countable families of solutions:

 {\rm (i)} A P-type family $\{F_k^P(z)>0, \, k \ge 0\}$ such that
  \beq
  \label{s23}
  F_k(0) \not = f_0 \quad \mbox{and} \quad F_k(z) = C_0 z^{-\frac
  {2d}{p-3}}(1+o(1)) \,\,\, \mbox{as} \,\,\, z \to + \infty,
   \eeq
   where $C_k>0$ is a constant, and each $F_k(z)$ has precisely
   $k+1$ intersections with the constant solution $f_0$; and

   {\rm (ii)} A Q-type family $\{F_k^Q(z), \, k \ge 0\}$, where each
   $F_k^Q(z)=f_0$ on some interval $z \in [0,a_k]$, $\{a_k>0\}$ is
    strictly monotone decreasing, with the
   following behaviour at the interface:
\beq
  \label{ff11N3}
 \mbox{$
F_k(z)= f_0 +(-1)^k \sqrt{\frac 89 \, \b a_k^3} \, (z-a_k)^{\frac
32}_+(1+o(1)) \quad \mbox{as} \quad z \to a_k,
 $}
 \eeq
 and has precisely $k$ intersections with $f_0$ for $z > a_k$ and
 has the asymptotic behaviour $(\ref{s23})$.


  \end{proposition}

\noindent For the main concepts of the proof, see  \cite{GPos, BuGa}.

\vspace{0.1in}

\begin{figure}
\centering
\subfigure[$k=0$]{
\includegraphics[scale=0.52]{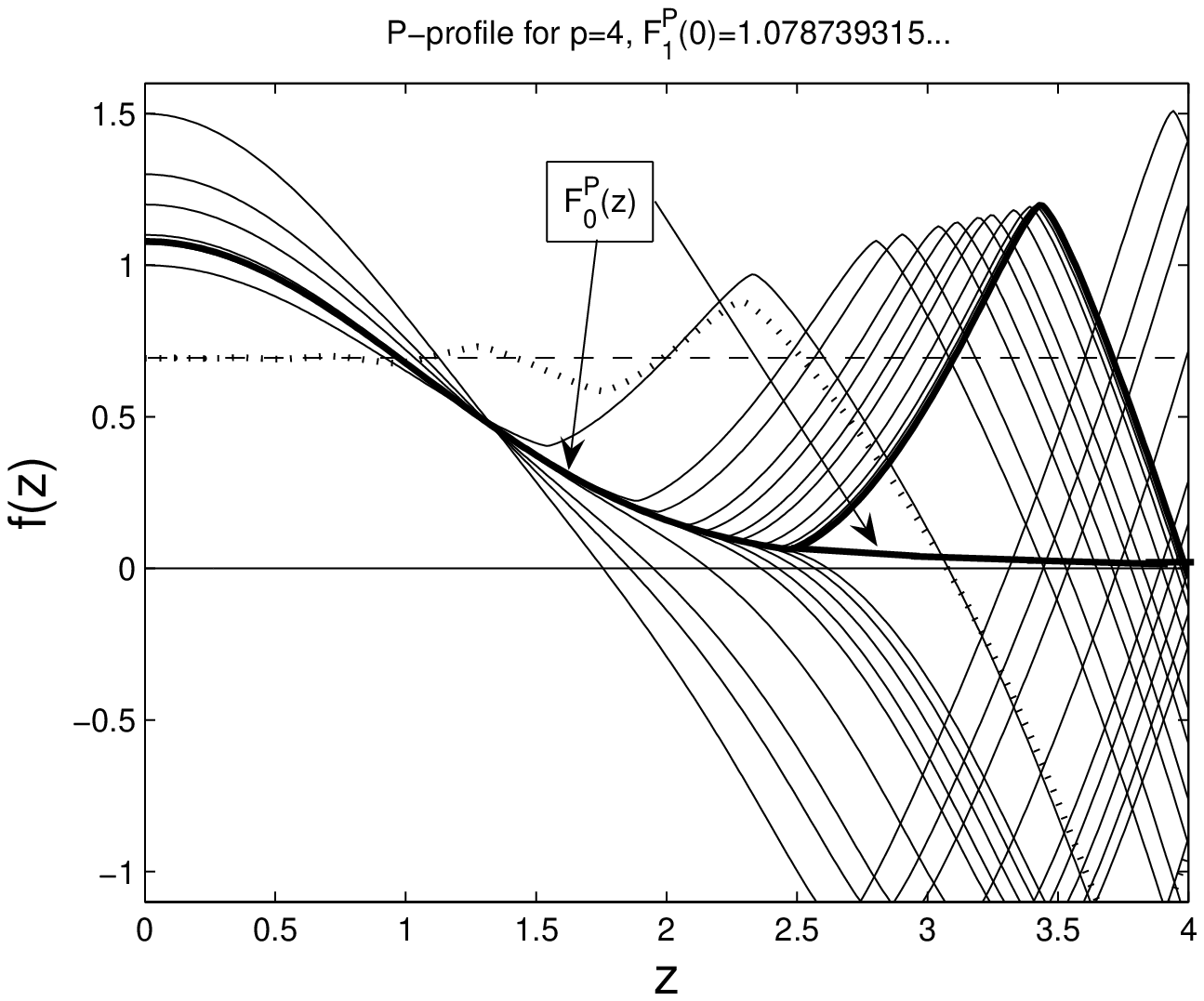} 
}
 \subfigure[$k=1$]{
\includegraphics[scale=0.52]{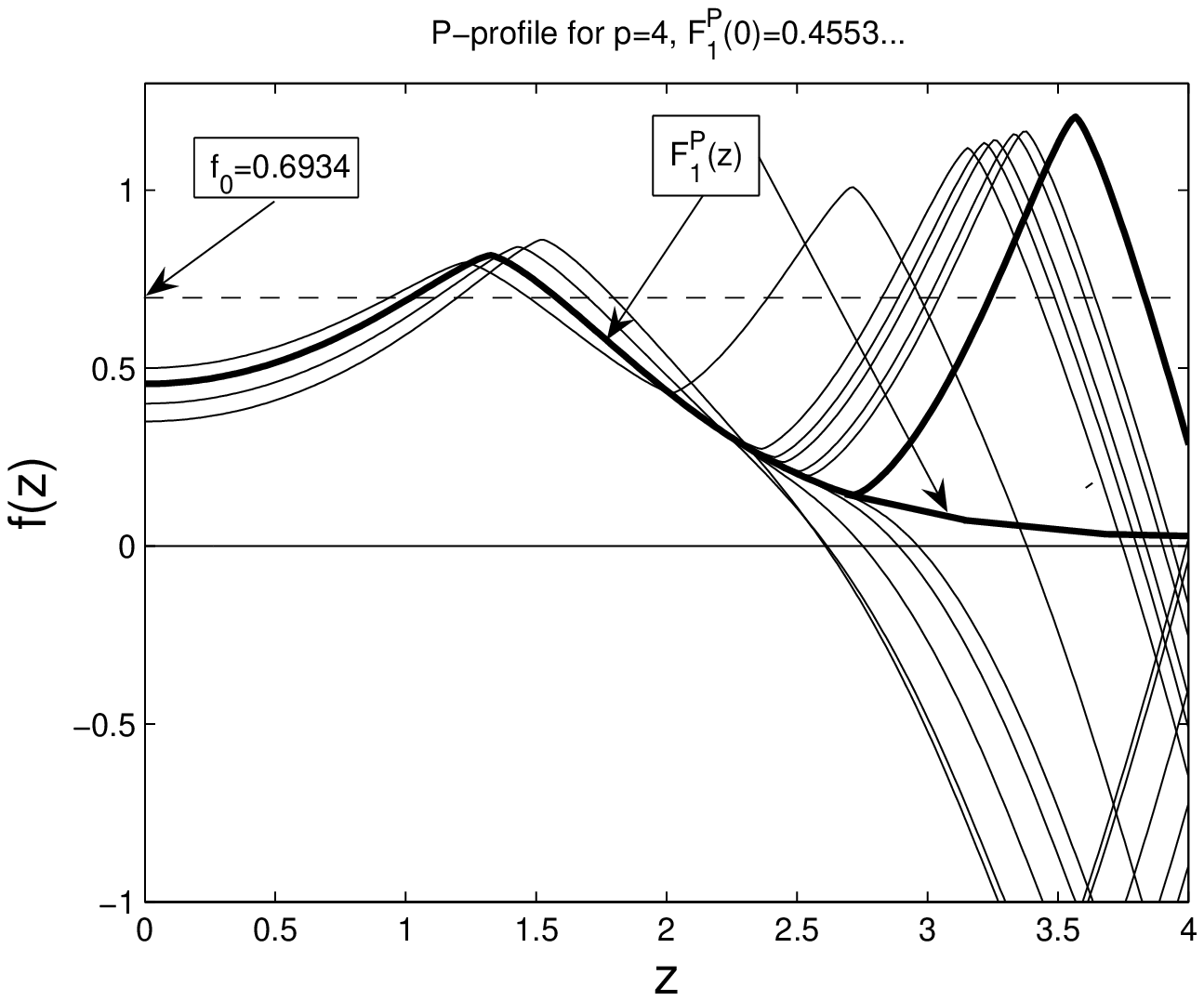} 
}
   \vskip -.3cm
    \caption{\small P-type profiles of the ODE (\ref{M6})
    for $d=3$, $p=4$.}
   \vskip -.3cm
 \label{F23}
\end{figure}

\noindent In Figure \ref{F23} (a,b), we show the similarity profiles $F_0^P(z)$
$F_1^P(z)$ for
$d=3, p=4$. Construction of further P-type profiles is similar
 and the following holds:
  $$
  F_k^P(0) \to f_0 \quad \mbox{as} \quad k \to \infty,
  $$
and moreover the convergence is from above for even $k=0,2,4,...$,
and from below for odd $k=1,3,5,... \, .$
 Uniqueness of each $F_k$ with $k+1$ intersections with equilibrium $f_0$
  is a difficult open problem.
We claim that as $k \to \infty$, both families $\{F_k^P(z)\}$ and
$\{F_k^Q(z)\}$ satisfying
 $$
 F_k^P(0) \to f_0 \quad \mbox{and} \quad a_k^Q \to 0
  $$
converge to a unique S-type profile $ F_\infty^S(z)$ such that
  \beq
  \label{fs1}
 F_\infty^S(0)=f_0 \,\,\, \mbox{and has infinitely many intersections with $f_0$ for small $z>0$.}
  \eeq

\vspace{0.1in}

\begin{figure}
\centering
\subfigure[$k=0$]{
\includegraphics[scale=0.52]{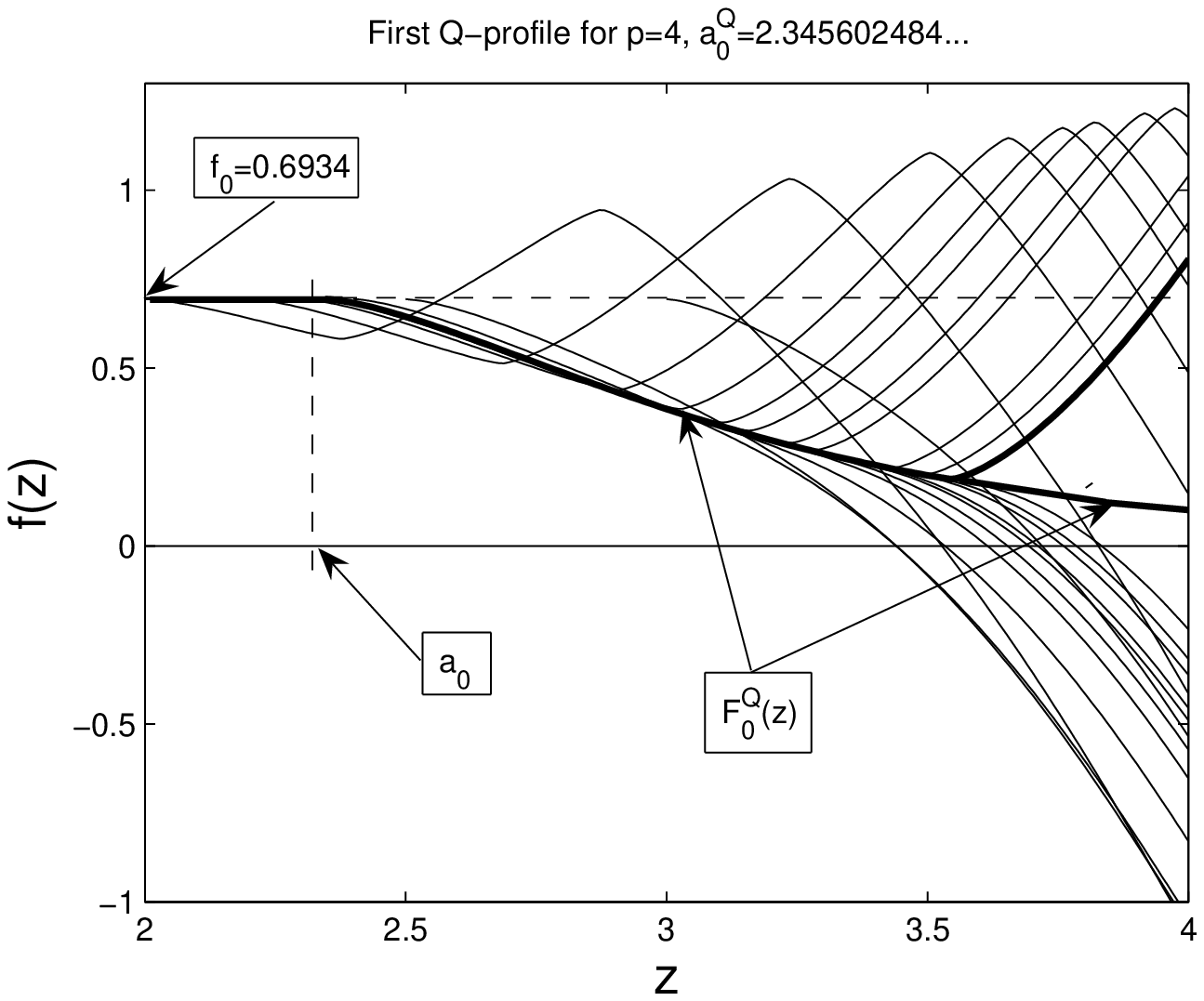} 
}
 \subfigure[$k=1$]{
\includegraphics[scale=0.52]{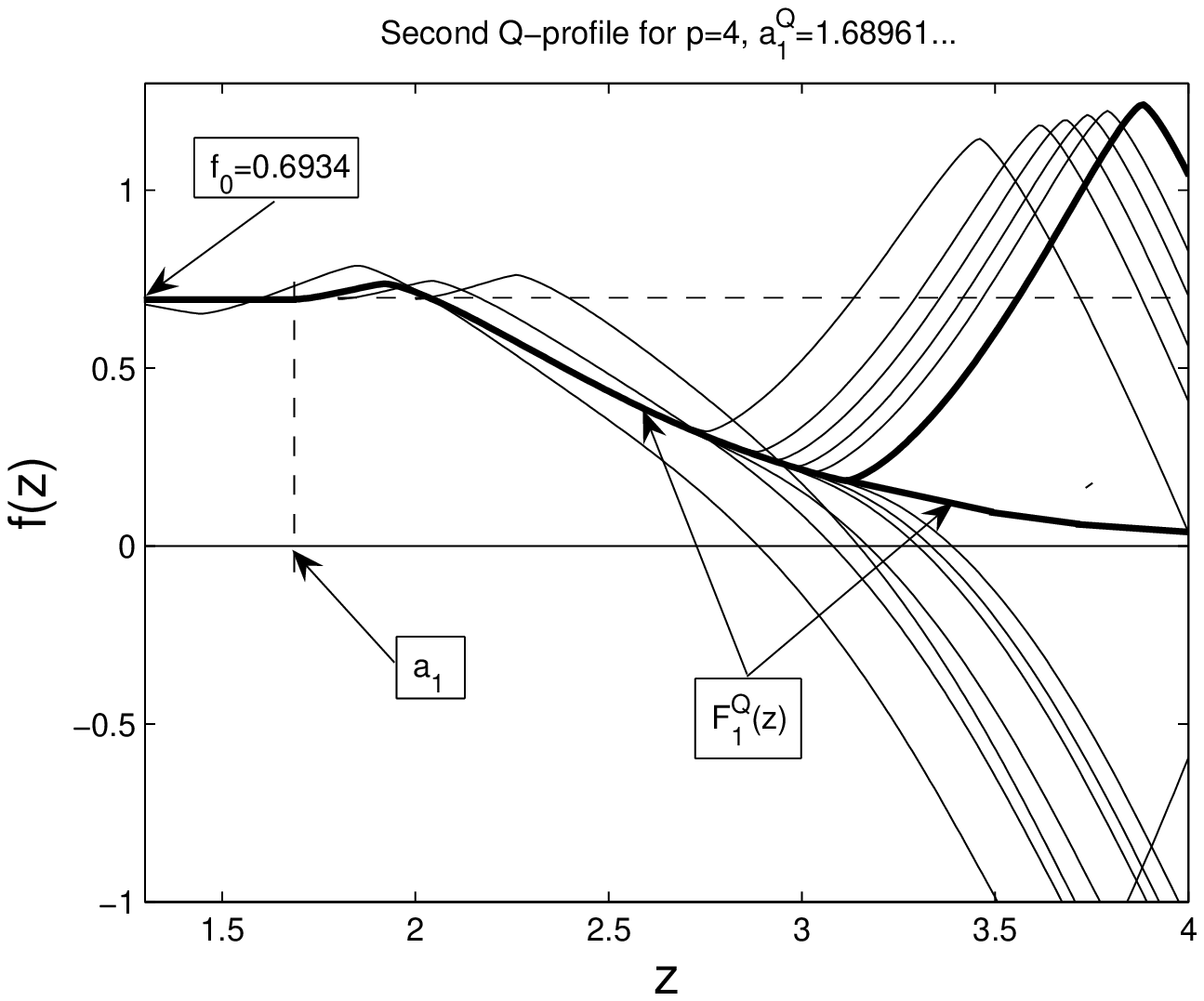} 
}
   \vskip -.3cm
    \caption{\small Q-type profiles of the ODE (\ref{M6})
    for $p=4$.}
   \vskip -.3cm
 \label{FQQ1}
\end{figure}

\noindent The first two similarity profiles of Q-type with $d=3$ and $p=4$,
with a flat centre part, are shown in Figure \ref{FQQ1}.

\vspace{0.2in}

\noindent Passing to the limit $t \to T^-$ in (\ref{M5}), where $z \to
 +\infty$, for any fixed $r>0$, we find from (\ref{s2}) the
 following final-time profile for both P and Q patterns:
  \beq
  \label{ft1}
  u_S(r,T^-)= C_0 r^{-\frac 2{p-1}} < \infty \quad \mbox{for any}
  \quad r>0,
  \eeq
  with $C_0>0$ fixed in (\ref{C01}).
  This implies  single point blow-up at the origin $r=0$ only
in both cases.

\vspace{0.1in}

\noindent We now make a {\em numerical study} of the stability of
these solutions with single-point blow-up. Using a similar method to that described in the
previous section (including taking Neumann boundary conditions
with $L = 8$) we can study the nature of the time dependent
blow-up solutions in this case. Usually when studying single-point
blow-up the narrowing of the solution peak as the reaction time $T$ is
approached, requires the use of a spatially adaptive mesh
\cite{BHR}. However in the M-A systems, as the blow-us the width
of the solution peak scales relatively slowly (as $(T-t)^{1/8}$
when $d=2$) it is not necessary to use an adaptive method to study
the nature of the blow up solutions in this case provided that the
spatial grid is fine enough. Taking $d=2$ and $p = 3$ and initial
data $u(r,t) = 10 {\mathrm e}^{-r^{2}}$ we find that the solution
blows up in a time $T = 0.028476$. Indeed, plotting $1/u^2(0,t)$
as a function of $t$ we find a close to linear solution, which has
a best fit with the equation $u(0,t) = 0.975/(T-t)^{1/2}$ with $T$
as above. Plotting the rescaled solution $(T-t)^{1/2}u(r,t)$ as a
function of the similarity variable $z = r/(T-t)^{1/8}$ we find
close agreement to the similarity solution constructed above. This
strongly implies that the monotone P-type similarity solution is
stable in the rescaled variables. A similar result is observed for
calculations when $d=3$. However, the Q-type and non-monotone
P-type blow-up profiles all appear from these calculations to be unstable.

\begin{figure}
\centering
\includegraphics[scale=0.4]{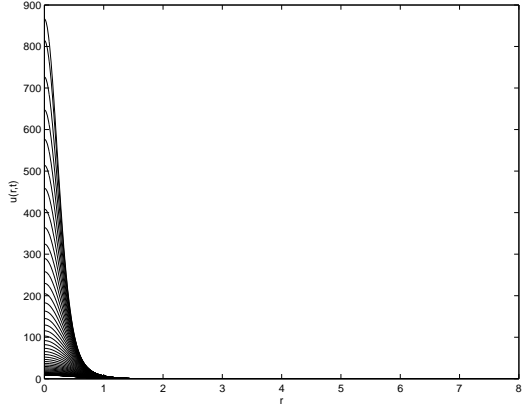}
\caption{\small The solution $u(r,t)$ when $d=2$ and $p=3$ as it evolves
toward a P-type singular solution at time $t = 0.028476$}
 \label{F2c}
\end{figure}

\begin{figure}
\centering
\includegraphics[scale=0.4]{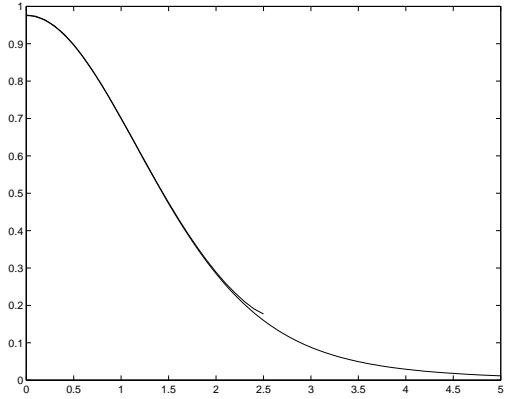}
\caption{\small The rescaled solution
$(T-t)^{1/2}u(r,t)$ plotted as a function
of the similarity variable $y = r/(T-t)^{1/8}$ compared
with the similarity profile of the P-type solution $F_0$
 obtained by shooting.}
 \label{F2cc}
\end{figure}


\ssk


\subsection{S-type periodic solutions}
\ssk

 \noindent The main part in the existence analysis for
(\ref{M6}), following the methods described in \cite{GPos, BuGa}, relies on constructing
two solutions to the initial value problem with {\em different oscillatory structure}, and
then deducing the existence of an intermediate solution with the correct properties at infinity
by  applying continuity arguments. To construct such solutions we must determine
the oscillatory properties of the solutions $f(z)$
about the positive equilibrium $f_0$. To study these we consider
solutions of the form
 \beq
 \label{l1}
 f(z)= f_0 + Y(z), \quad \mbox{where $Y$ is small and solves the reduced equation}
  \eeq
   \beq
   \label{l2}
    \mbox{$
   \frac{1}{z^{d-1}} \, |Y'|^{d-1} Y'' - \b Y' z+Y=0,
    $}
    \eeq
    which remains nonlinear. In order to study oscillations of
    $Y(z)$ about 0, we can exploit the
     scaling invariance of (\ref{l2}) and introduce the {\em oscillatory
     component}
  $\varphi$   as follows:
  \beq
  \label{l3}
   \mbox{$
  Y(z) = z^{\alpha} \varphi(s), \quad s= \ln z, \quad \mbox{where}
  \quad  \alpha = \frac{2d}{d-1}.
   $}
   \eeq
 Substituting (\ref{l3}) into (\ref{l2}) in (for example) the case of $d=2$
yields the following {\em autonomous ODE}:
  \beq
  \label{l4}
   \mbox{$
 |\vp'+4 \vp|(\vp''+7 \vp'+12 \vp)- \b \vp'+ \frac 1{p-1} \, \vp=0.
  $}
  \eeq
Setting $\vp'=\psi(\vp)$ yields the first-order ODE system
 \beq
 \label{l5}
\psi \frac{{\mathrm d} \psi}{{\mathrm d}\vp}=
\frac{\b \psi -\frac 1{p-1} \, \vp}{|\psi + 4\vp|} -7 \psi-12\vp,
\eeq
which we can study by using phase-plane analysis, and particular identify a limit cycle.

\begin{proposition}
 \label{Pr.lc}
 The ODE $(\ref{l5})$ admits a stable limit cycle on the
 $\{\vp,\psi\}$-plane, which generates a periodic solution $\varphi_*(s)$ of $(\ref{l4})$.
  \end{proposition}

\noi{\em Proof.} Equations (\ref{l5}) from PME and $p$-Laplacian
theory with limit cycles have been studied since 1980's; see
\cite{GPos} and extra references and related results in
\cite{BuGa}. The limit cycle exists for all $p>d$. $\qed$

\ssk

\noindent As an immediate consequence,
the gradient-dependent ODE (\ref{M6})
also admits an S-type solution $F^S(z)$ with infinitely many oscillations about
 the equilibrium $f_0$ near the origin.
 Since the behaviour (\ref{l1}), (\ref{l3}) violates the monotonicity,
this S blow-up profile is not associated with the original M-A
equation when $d=2$ but does correspond to a possible solution
when $d=3$. \noindent The  existence of $F^S$ then follows by
shooting from the origin using the
  1D asymptotic bundle constructed using (\ref{l1}) and the
   representation (\ref{l3}) with the periodic solution $\varphi_*(s)$, i.e.,
  \beq
  \label{ff11}
  f(z) = f_0 + z^{\frac{2d}{d-1}} \varphi_*(s_0+ \ln z) +... \, .
 \eeq
Here the shift $s_0$ is the only shooting parameter. Using this parameter,
 it follows from the  oscillatory structure of the expansion (\ref{ff1})
that there
  exists an $s_0$ such that
 $f=F^S$ has the power decay at infinity given by (\ref{s2}) with a positive constant $C_0$.
  We do not have clear evidence for the uniqueness of such $F^S(z)>0$
and it appears from the numerical calculations that the associated
self-similar blow-up profiles are unstable.


\subsection{Global blow-up for $p \in(1,d)$}

\noindent Many of the mathematical features of the analysis of the similarity
blow-up structures in this case are the same as for $p>d$, though the evolution
properties are quite different. It follows from (\ref{M5}) that $p
\in(1,d)$ corresponds to the {\em global blow-up}, where
 \beq
 \label{l6}
 u_S(x,t) \to \infty \quad \mbox{as \, $ t \to T^-$
 \, uniformly on bounded intervals in $\re_+$}.
  \eeq
 As in \cite[pp.~183-189]{SGKM}, we obtain existence
  of a similarity profile. Since the PDE (\ref{2.1}) is non-autonomous in space,
  uniqueness remains an open problem,
  since  the geometric Sturmian approach to uniqueness (see \cite{CGF}
  for main results and references) does not
   apply.

\begin{proposition}
 \label{Pr.HS}
 For any $p \in (1,2)$,
  problem $(\ref{M6})$ admits a  compactly supported
 solution $F_0(y)$.
 \end{proposition}

\begin{figure}
\centering
\includegraphics[scale=0.7]{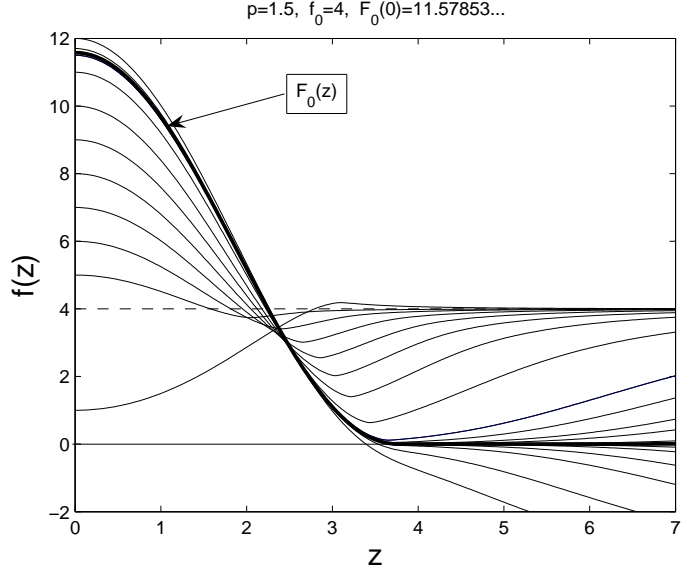}
\caption{\small Shooting compactly supported similarity profile
$F_0(z)$ (the bold line) for $p=1.5$ in (\ref{M6});
$F_0(0)=11.5785...$\,.}
 \label{F4}
\end{figure}

\noindent The numerical shooting construction of $F_0$ is shown in Figure
\ref{F4} by the bold line for $d=2$ and $p= \frac 32$.

 \section{Examples of blow-up in a non-radial geometry in $\re^2$}
\label{SectNR}

\noindent It is immediate that the M-A equation (\ref{M1}) admits
non-radially symmetric blow-up solutions which do not become more
symmetric as $t \to T^-$. These solutions can be obtained directly from a
non-symmetric transformation of the radially symmetric blow-up
solutions described in the previous section and we describe below. However, it is
unclear from the analysis whether these solutions are stable or
not. In this section we consider these solutions and make some
numerical computations to infer their stability.


\subsection{Regional blow-up for $p=d=2$: the existence of non-radially symmetric blow-up solutions}

 For $p=d=2$, (\ref{M1}) is
 \beq
 \label{M1NN}
  u_t = - | D^2 u| + u^2 \quad \mbox{in} \quad \re^2 \times \re_+.
   \eeq
This partial differential equation admits self-similar blow-up
solutions  (\ref{M5}), i.e.,
 \beq
 \label{R0}
  \mbox{$
 u_S(x,y,t) = \frac 1{T-t} \, f(x,y).
  $}
  \eeq
Here the function $f(x,y)$ now satisfies the following ``elliptic"
M-A equation, with decay to zero at infinity:
 \beq
 \label{R1}
 {\bf A}(f) \equiv -|D^2 f|+ f^2 - f=0 \quad \mbox{in} \quad \re^2.
   \eeq
The radial compactly supported solution $F_0(r)$ described in
Proposition \ref{Pr.1} also solves (\ref{R1}). The existence of
non-radial smooth solutions of (\ref{R1}) now follows from an
invariant group of scalings. Indeed, if $f(x,y)$ is a solution of
(\ref{R1}) then so is the function
 \beq
 \label{R66}
  \mbox{$
 f_a(x,y)= f\big(\frac x a, a y\big) \quad \mbox{for any\,
 $a>0$}.
  $}
  \eeq
Indeed any rotation of this function is also a solution.
 Therefore, by taking the radial profile $F_0(y)$ supported in $[0,L_S]$ as a solution of
 (\ref{M7}), we can obtain the non-radially symmetric blow-up solution
  \beq
  \label{R67}
   \mbox{$
  u_{\rm S}(x,y,t) = \frac 1{T-t} \, F_0 \big(\sqrt{(\frac x a)^2 +
  (a y)^2}\,\big),
   $}
    \eeq
 This blows up in the ellipsoidal {\em localization domain} given by
 \beq
 \label{R68}
 E_a = \big\{ (x,y):  \quad \mbox{$
\big(\frac x a \big)^2 +
  (ay)^2 < L_S^2 \quad (a \not = 1) \big\}.
   $}
   \eeq
   Note that its area does not depend on $a$:
 $$
 {\rm meas} \, E_a = \pi L_S^2 \quad \mbox{for any} \,\,\, a>0.
  $$
Thus, M-A equations such as (\ref{M1NN}) do not support the
(unconditional) symmetrization phenomena, found for many
semilinear and quasilinear parabolic equations; see
\cite[p.~50]{AMGV} for references and basic results. In classic
parabolic theory, results on symmetrization are well known and
  are connected with  the moving plane method
and Aleksandrov's Reflection Principle; see key references in
\cite[Ch.~9]{GilTr} and \cite[p.~51]{AMGV}. However, all these
approaches are based on the Maximum Principle that fails for M-A
parabolic flows like
  (\ref{M1NN}).
  We conjecture that (\ref{R1}) does not admit other non-symmetric
solutions but have no proof of this result.

\subsection{On the linearized operator}

 \noindent Checking the stability properties of the non-radial solutions of
 (\ref{R1}), one can easily derive the linearized operator about
 the radial state $F_0(r)$
  \beq
  \label{L1}
    \begin{matrix}
   {\bf A}'(F_0)Y= \big[F_0'' \frac{y^2}{r^2}+ F_0'\big(\frac
   1r-\frac{y^2}{r^3}\big)\big]Y_{xx}+
   \big[F_0'' \frac{x^2}{r^2}+ F_0'\big(\frac
   1r-\frac{x^2}{r^3}\big)\big]Y_{yy} \ssk\ssk  \\
   -2\big(F_0''- \frac 1r\,F_0'\big) \frac{xy}{r^2}\, Y_{xy}-Y.
   \end{matrix}
    \eeq
 Moreover, since ${\bf A}$ is potential in $L^2$, this Frechet
 derivative is symmetric, so there is a hope to get a ``proper"
 self-adjoint extension of the linear operator (\ref{L1}).
 Unfortunately, we should recall that (\ref{L1}) cannot be treated
 as  elliptic in the domain of convexity of $F_0$.
 In addition, since $F_0(r)$ is compactly supported, the operator
 (\ref{L1}) have  singular
 coefficients at $r=L_S$ and will be inevitably defined in
 a complicated domain with possibly
  singular weights,
 which makes rather obscure using such
 operators in studying the angular stability or unstability of the
 radial profile $F_0(r)$. In any case, it is convenient to note
 that, due to the symmetry (\ref{R66}), the stability is {\em neutral},
 i.e.,
  \beq
  \label{L2}
   \mbox{$
  \exists \,\, \hat \l=0, \quad \mbox{with the non-radial eigenfunction}
  \quad
  \hat \psi_0=  \frac{\mathrm d}{{\mathrm d}a}\, f_a\big|_{a=1}=-
  \frac{x^2-y^2}{r}\, F'_0(r)
  $}
  \eeq
   (we naturally assume that the eigenfunction belongs to  the domain of the
   self-adjoint extension). In other words, the
   stability/unstability will depend on an appropriate and delicate centre
   manifold behaviour. More precisely, if we perform the linearization
    $f(\t)=F_0+Y(\t)$ of the
   corresponding to (\ref{R1}) non-stationary flow
    \beq
    \label{L3}
     f_\t= {\bf A}(f)
      \quad \Longrightarrow \quad
      Y_\t={\bf A}'(F_0)Y-|D^2 Y|+Y^2.
       \eeq
  Then the corresponding formal centre subspace behaviour\footnote{Existence of a centre manifold
  by standard invariant manifold theory \cite{Lun} is a very difficult
 open problem, which seems hopeless.} by setting
  $$
  Y(\t)= a(\t) \hat \psi_0 + w^\bot, \quad w^\bot \bot \hat \psi_0 \quad \big(\|
w^\bot(\t)\| \ll
  |a(\t)|\,\,\,\mbox{for} \,\,\, \t \gg 1\big)
   $$
   leads, on projection (as in classic theory, we have to assume at the moment
    a certain completeneee/closure of
   the orthonormal eigenfunction subset, which are very much questionable problems), to the
equation
 \beq
 \label{L5}
  \dot a= \g_0 a^2 \quad \mbox{for $\t \gg 1$, where}
  \quad \g_0=\langle -|D^2 F_0|+F_0^2, \hat \psi_0 \rangle.
   \eeq
Therefore, for $\g_0>0$, stability/instrability of such flows
depend on how the sign of $\g_0$ is associated with the sign of
the expansion coefficient $a(\t)$.

More careful checking by using equation (\ref{R1}) and the
eigenfunction in (\ref{L2}) of changing sign shows that
 \beq
 \label{L6}
  \g_0= \langle F_0, \psi_0 \rangle=0,
   \eeq
 so that this centre subspace angular evolution according to (\ref{L5}) is
 formally absent at all. Of course, this is not surprising, since our
 functional setting assumes fixing the domain $\{r < L_S\}$ of
 definition of the functions involved, and this clearly prevents
 any angular evolution on the centre subspace that demands
 changing this domain.


In view of such a non-justifying formal linearized/invariant
manifold analysis, we will next rely on also rather delicate
numerical techniques to check angular stability of blow-up
similarity profiles and solutions.


 \subsection{Single point blow-up in non-radial geometry: similarity solutions for $p>d=2$}
\label{SectNR2}


\noindent We can use a similar method to study non-radially symmetric single-point
blow-up profiles. The self-similar solution
(\ref{M6}),
  \beq
 \label{uu0}
  \mbox{$
 u_S(x,y,t)= (T-t)^{-\frac 1{p-1}} f(\xi,\eta), \quad
 \xi=x/(T-t)^\b, \,\,\, \eta=y/(T-t)^\b, \quad \b=
 \frac{p-2}{4(p-1)},
  $}
  \eeq
 now leads to a more
complicated elliptic M-A equation
 \beq
 \label{uu1}
  \mbox{$
 -|D^2 f|- \b \n f \cdot \z - \frac 1{p-1} \, f + |f|^{p-1}f=0
  \quad \mbox{in} \quad \re^2 \quad \big(\z=(\xi,\eta)^T\big).
   $}
   \eeq

\noindent As before, the group of scalings (\ref{R66}) leaves
equation (\ref{uu1}) invariant. Therefore, taking the radial
solution $F_0(z)$ from Proposition \ref{PrLS}, we obtain the
family
 \beq
 \label{ff1}
  \mbox{$
 F_a(\xi,\eta) = F_0\big(\sqrt{(\frac \xi a)^2 + (a \eta)^2}\, \big)
  $}
  \eeq
  of non-radial solutions  of (\ref{uu1}) (together with
all rotations of these). However, this set of solutions may not
exhaust all the non-symmetric patterns. To see this, we consider
the linearization (\ref{l1}) about the constant equilibrium $f_0$
which leads to a nonlinear M-A elliptic problem:
 \beq
 \label{uu2}
-|D^2 f|- \b \n f \cdot y + f \equiv
-\big[f_{\xi\xi}f_{\eta\eta}-(f_{\xi\eta})^2\big]- \b (f_\xi \xi +
f_\eta \eta) + f =0
  \quad \mbox{in} \quad \re^2.
   \eeq
This fully nonlinear PDE does not admit separation of variables.
We can see from (\ref{uu1}) that if $f(z) \to 0$ as $z \to \infty$
  sufficiently fast, the far-field behaviour is
 governed by the linear terms,
 \beq
 \label{uu3}
  \mbox{$
 - \b
(f_\xi \xi + f_\eta \eta) - \frac 1{p-1} \, f +...=0.
 $}
  \eeq
  Solving this gives the following typical asymptotics (cf. (\ref{s2})):
   \beq
   \label{uu4}
    \mbox{$
   f(z)= C\big(\frac z{|z|}\big) |z|^{- \frac 4{p-2}} + ... \quad \big(z=(\xi,\eta)^T\big),
    $}
 \eeq
 where $C(\mu) > 0$ is an arbitrary smooth function on the unit circle
 $\{|\mu|=1\}$ in $\re^2$. The constant function $C_0(\mu) \equiv
 C_0>0$ gives the radially symmetric similarity profile as in Proposition
 \ref{PrLS}. Furthermore the $\pi$-periodic function $C(\mu)$
 given in the polar angle $\varphi$ by
  $$
  \mbox{$
  C_1(\varphi)= \frac 12\big(\frac 1{a^2}+a^2\big) + \frac 12
  \big(\frac
1{a^2}-a^2\big) \cos 2\varphi
  $}
  $$
 generates the ellipsoidal
 solutions (\ref{ff1}).
We conjecture that other solutions are possible
with $C_l(\varphi)$ having smaller periods
 $\frac {2 \pi}3$, $\frac \pi 2$,\,...\,.
However, at present, the existence of such is unknown.



\subsection{Numerical computations of the non-radially symmetric
time-dependent solutions}

\noindent We now consider a numerical computation of the blow-up
profiles when $\Omega = [0,1] \times [0,1]$ is the unit square (so
that $d=2$), and we took $p=3$. It is convenient in this
calculation to impose Dirichlet boundary conditions. The PDE is
solved by using a semi-discrete method for which the square is
divided into a uniform grid (typically a $100 \times 100$ mesh)
and the spatial Monge--Amp\`ere operator discretised in space by
using a second-order 9 point stencil. The resulting time dependent
ODE system is very nonlinear and an implicit solver is very
inefficient. Accordingly it was solved using an explicit, adaptive
Runge--Kutta method with a small tolerance. The discretisation in
space leads to certain chequer-board type
 instabilities\footnote{We recall that the M-A flow under consideration is supposed to have
 some natural instabilities in the areas, where the concavity of solutions
 is violated; these questions will be discussed.} and these
are filtered out at each stage by using a suitable averaging
spatial filter applied to the ODE system. For initial data
satisfying the Dirichlet condition we took
  $$
   \mbox{$
  u(0,t) = 10^4 {\mathrm e}^{-4
r^{2}} \sin(\pi x) \sin( \pi y), \quad \mbox{where} \quad r^{2} =
a^2(\hat{x}-0.4)^{2} + \frac 1{a^2}\, (\hat{y}-0.6)^{2}, $}
 $$
for which $a = 2$ and $\hat{x},\,\,\hat{y}$ were a set of
coordinates rotated at an angle of $\frac  \pi 4$. This data were
chosen to have an elliptical set of contours close to its peak.

\begin{figure}
\centering
\subfigure[Contours]{
\includegraphics[scale=0.52]{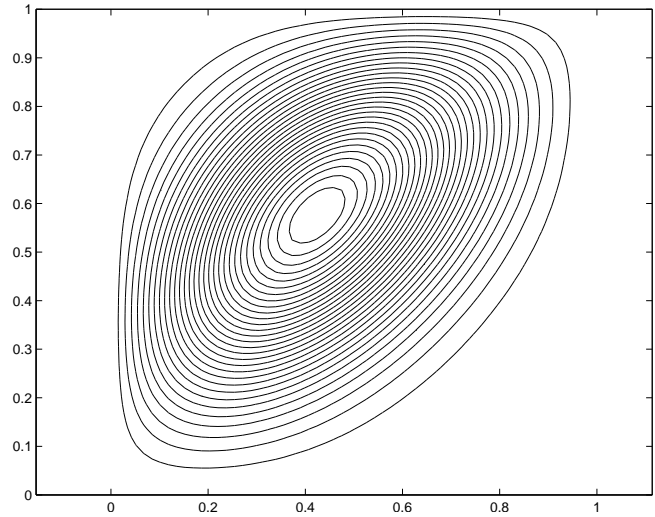} 
}
 \subfigure[Profile]{
\includegraphics[scale=0.52]{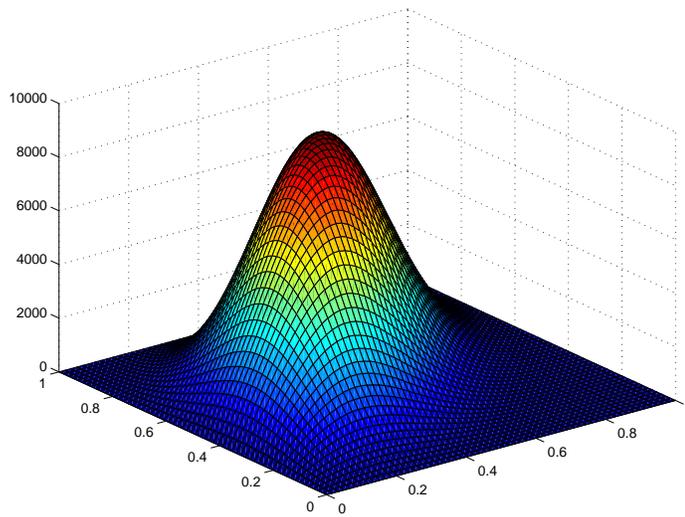} 
}
   \vskip -.3cm
    \caption{\small Initial solution profile and contours for $d=2$, $p=4$.}
   \vskip -.3cm
 \label{Fc1}
\end{figure}

\begin{figure}
\centering
\subfigure[Contours]{
\includegraphics[scale=0.52]{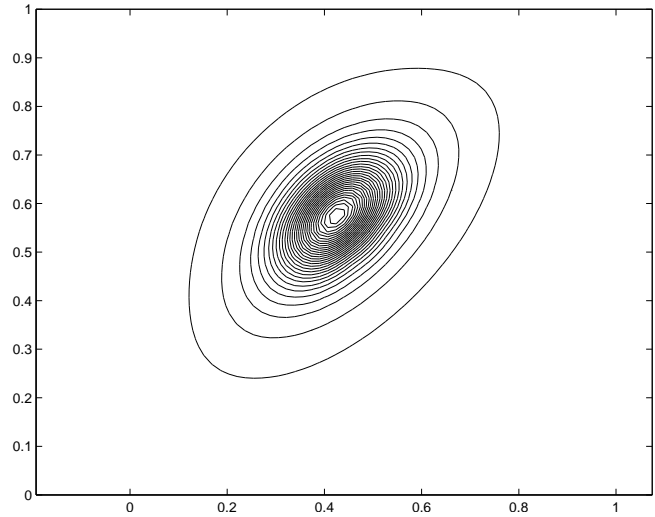} 
}
 \subfigure[Profile]{
\includegraphics[scale=0.52]{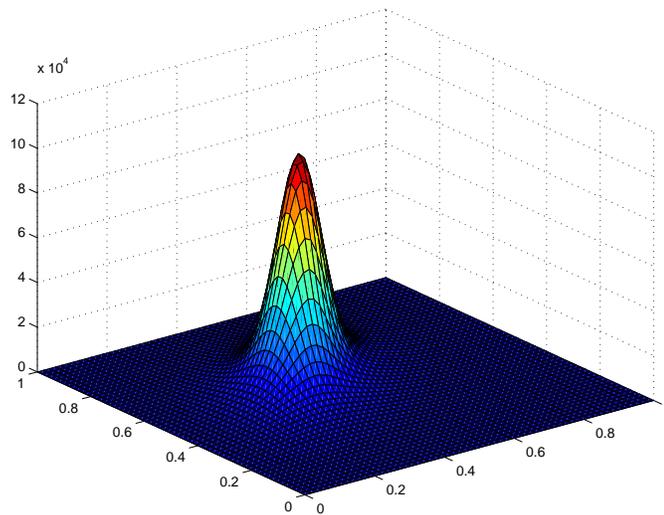} 
}
   \vskip -.3cm
    \caption{\small Profile and contours of the solution for $ d=2, p=4$ closer
to the blow-up time.}
   \vskip -.3cm
 \label{Fc2}
\end{figure}

\begin{figure}
\centering
\includegraphics[scale=0.6]{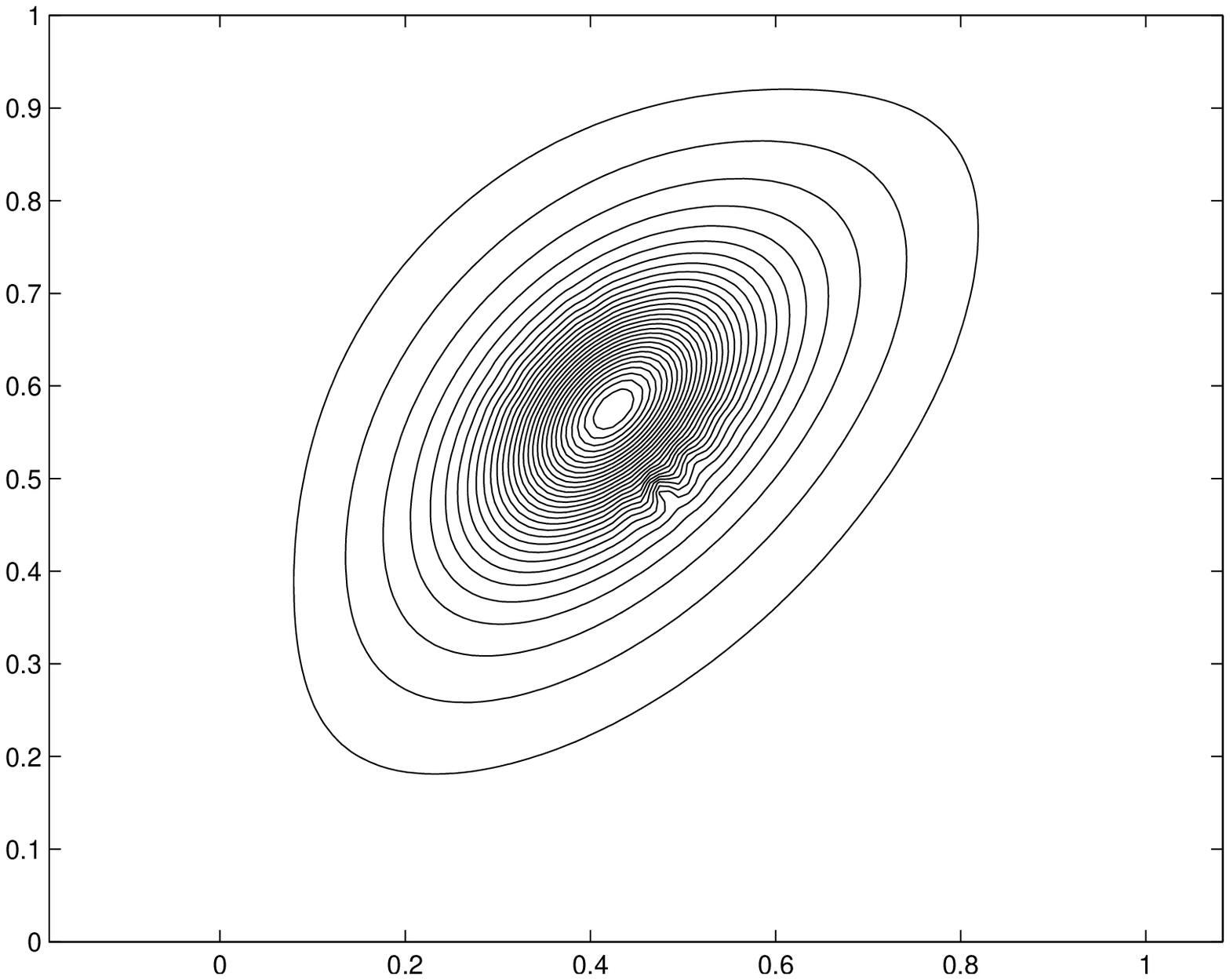}
\caption{\small The solution at a slightly later time showing a
possible instability at as point where the profile loses
concavity.}
 \label{Fc3}
\end{figure}

This system led to blow-up of the discrete system in a computed finite time
$T \approx 1.17 \times 10^{-4}$. (Note that the computed blow-up
time decreases when the spatial mesh is refined).
In Figure~\ref{Fc1} (a,b) we present the initial solution and its contours
and in Figure~\ref{Fc2} (a,b) a solution much closer to the blow-up time
(so that it is approximately 10 times larger than the initial
profile).
Note that the elliptical form of the contours has been preserved
during the evolution giving some evidence for the stability of
the elliptical blow-up patterns. We also give in Figure~\ref{Fc3}
a plot of the solution a slightly later time. Note in this case
evidence for an instability close to a point where the solution profile
loses convexity. It is not clear at present whether this is a numerical or
a true instability. Certainly all of the numerical methods used had
extreme difficulty in computing a significant way into the blow-up
evolution.

\section{On conservation of mass, comparison, existence-uniqueness, etc.}
 \label{SectEx}

Here
 we briefly discuss some related analytical theoretical aspects  M-A
equations such as (\ref{M1NN}), which appeared before and have
been already slightly discussed some times. We skip the quadratic
reaction term $u^2$, and concentrate on the fully nonlinear M-A
operator, for which we consider the Cauchy problem
 \beq
 \label{R4}
 u_t= -|D^2 u| \equiv -u_{xx}u_{yy} + (u_{xy})^2 \quad \mbox{in} \quad \re^2 \times \re_+, \quad
 u(x,0) = u_0(x) \ge 0 \quad \mbox{in} \quad \re^2,
  \eeq
with sufficiently smooth compactly supported initial data
satisfying some extra necessary conditions (e.g., of having
``dominated concavity").
 Since
compactly supported solutions are supposed to be involved, we
cannot use the advanced techniques of classic nowadays theory of
convex (concave) solutions; see \cite{GilTr, Iv94}, etc.
 Notice that, even for convex solutions, local regularity theory
 for the basic M-A equation
  $$
   {\rm det} \, D^2 u=f(x)>0 \quad \mbox{in a convex bounded domain}
   \quad \O \subset \ren
 $$
 is rather involved with a number of open questions; see
 \cite{Huang06} as a guide.
 On the other hand,  even for 2D stationary M-A  equations of
 changing convexity-concavity (see (\ref{Kh1})), there are counterexamples
 concerning local solvability and regularity, \cite{Khuri07}.
 These questions of M-A theory have not been developed in detail
 still\footnote{``Yet, it is remarkable that the basic question of whether
 there exist any examples of local
  nonsolvability, has remained open for this well-studied class of
  equations",
  \cite[p.~665]{Khuri07}.}.
  The study of
 finite regularity solutions of the simplest degenerate (at 0) M-A equation in
 $\re^2$ with radial homogeneous $f(x)$,
   \beq
   \label{DDDD1}
{\rm det} \, D^2 u=|x|^\a \quad \mbox{in}\quad  B_1
  \eeq
has some surprises \cite{Dask07}; e.g. for $\a>0$, there exist a
radial and a non-radial $C^{2,\d}$ solutions, for $\a \in (-2,0)$
the only radial does (this is about a delicate study of a single
point  blow-up singularity for (\ref{DDDD1})); see also
\cite{Rios08} for regularity of radial solutions for (\ref{DDDD1})
with the right-hand side $f(\frac 12\, |x|^2,u, \frac 12 \, |\n
u|^2)$. Equation (\ref{DDDD1}) has the origin in Weyl's classic
problem (1916).

   We do not plan to discuss such delicate questions somehow seriously here
   (especially for our problems having so strong degeneracy and even changing of type),
 and restrict to first auxiliary  aspects of such singularity
 phenomena for (\ref{R4}).

\subsection{Source-type similarity solutions}

These are easy to construct for the radial equation
 \beq
 \label{81}
  \mbox{$
 u_t = - \frac 1r \, u_r u_{rr}, \quad \mbox{so that}
 $}
 \eeq
 \beq
 \label{82}
  \mbox{$
 u_*(r,t)= t^{-\frac 13} F(y), \quad y=\frac r{t^{1/6}}, \quad
 \mbox{where} \quad F(y)= \frac 1{48} \,\big(d^2-y^2\big)^2_+, \,\,\, d>0.
  $}
  \eeq

 \subsection{Scaling group: non-symmetric solutions}

 Let us mention that, obviously, (\ref{R4}) admits a variety of
 non-radial solutions. Indeed, the equation is invariant under the
 following group of scaling transformations:
  \beq
  \label{sc1}
   \mbox{$
   u(x,y,t) \mapsto \frac {a^2 b^2}c \, u\big( \frac x a, \frac y
   b, \frac t c \big), \quad a, \, b, \, c \not = 0.
    $}
    \eeq
Therefore, (\ref{R4}) does not support the symmetrization
phenomena that are typical for many nonlinear parabolic PDEs.

Using also the time-translation, we obtain from (\ref{82}) the
following 4-parametric family of exact solutions:
 \beq
 \label{zz1}
  \mbox{$
 u_*(x,y,t)=\frac{a^2 b^2}{c^{2/3}} \, (T+t)^{-\frac 13}
 \big[d^2 - c^{\frac 13} (T+t)^{-\frac 13}
 \big( \frac{x^2}{a^2} + \frac{y^2}{b^2}\big) \big]^2_+.
 $}
  \eeq

\subsection{No order-preserving semigroup in non-radial geometry}

We recall that in the radial geometry, the semigroup for the
parabolic equation (\ref{81}) is obviously order-preserving. It
turns out that in the non-radial setting, this is not the case.

\begin{proposition}
  \label{Pr.No}
  In general, sufficiently smooth solutions of  $(\ref{R4})$ do not obey
  comparison.
   \end{proposition}

\noi{\em Proof.} We take two exact solutions from (\ref{zz1}):
$u(x,y,t)$ with $a=b=c=T=1$ and the general solution $u_*(x,y,t)$,
and show that the usual comparison is violated in this family of
non-radial solutions. Comparing positions of the interfaces at the
$x$- and $y$-axes and the maximum values at the origin yields for
initial data at $t=0$ that
 \beq
 \label{y1}
  \mbox{$
 u_*(x,y,0) \le u(x,y,1)  \quad \mbox{if}
 \quad  \frac {a d}{c^{1/6}} T^{\frac 16}<1, \,\,
\frac {b d}{c^{1/6}} T^{\frac 16}<1, \,\,\frac {a^2 b^2
d^4}{c^{2/3}} T^{-\frac 13}<1.
 $}
 \eeq
On the other hand,  the comparison is violated for large $t \gg1$
if
 \beq
 \label{y2}
 \mbox{$
 \frac{ad}{c^{1/6}} \, (T+t)^{\frac 16} > t^{\frac 16}, \quad
 \mbox{i.e.,} \quad \frac{ad}{c^{1/6}} > 1.
 $}
  \eeq
 It is easy to see that the system of four algebraic inequalities
in (\ref{y1}) and (\ref{y2}) has, e.g., the following solution:
 $$
  \mbox{$
 a=1, \quad b= \frac 18, \quad c=2, \quad d=2, \quad T \in
 \big(\frac 1{2^8}, \frac 1{2^5}\big). $\qed$
  $}
  $$

 \subsection{Towards well-posedness}

 We present the following very formal speculations, which nevertheless  show that
the non-fully concave M-A flow (\ref{R4}) is better well-posed as
it can be expected. Actually, this was observed in a number of
numerical experiments discussed above. Assume that, due to an
essentially deformed spatial shape of the solution (say, by means
of choosing special ``ellipsoidal" initial data), we consider the
{\em unstable area} that is characterized as follows: $|u_{xy}|^2
\ll |u_{xx}u_{yy}|$ and, e.g., $u_{yy} \ge c_0>0$, i.e., the flow
 \beq
 \label{fl1}
 u_t =- u_{yy} u_{xx}+...\, , \quad u(x,y,t)>0 \quad\big(\mbox{cf.} \quad u_t=-c_0
 u_{xx}+...\big)
  \eeq
  is backward parabolic with respect to the spatial variable $x$. Then,
   let us assume that the positive solution $u(x,y,t)$ is going to
   produce  a blow-up singularity in finite time as $t \to T^-$, and, say, let it
   be a Dirac's delta of a positive measure (blow-up of derivatives is a different matter, which
   sometimes can be also  treated). Of course, we do not
   mean precisely that in this fully nonlinear equation, but can
   expect that a certain such tendency as $t$ moves towards $T$ can be
   observed, as the linear PDE in the braces in (\ref{fl1}) suggests.
   Hence, if such a tendency of approaching a $\sim \d(x-x_0)$ in $x$ is
   observed, then, obviously, at this unstable subset
 \beq
 \label{fl2}
 u_{xx} \le -c_0<0 \quad \Longrightarrow
 \quad  u_t =(-u_{xx}) u_{yy} +...\,\,\, \mbox{becomes well-posed
 parabolic in $y$}
  \eeq
 (here we again assume that $u_{xy}$ does not play a role at this
 stage).
 In other words, such a simple localized pointwise singularity
 $\sim \d({\bf x}-{\bf x}_0)$ in both variables $x$ and $y$ is
 unlikely.
 This means that the PDE (\ref{R4}) can exhibit a certain
 ``self-regularization" even in the case of not fully concave data.
 We are not aware of any rigorous mathematical justification of
 such a phenomenon, and will continue to discuss this subject
 below using other arguments.

\subsection{$\e$-regularization: on formal extended semigroup}

We propose to construct a unique {\em proper} solution of
(\ref{R4}) as a limit of a family of smooth regularized solutions
$\{u_\e\}$ satisfying the regularized fourth-order uniformly
parabolic equation
 \beq
 \label{R5}
 u_\e: \quad u_t=- \e \D^2 u -|D^2 u| \quad \mbox{in} \quad \re^2 \times \re_+,
  \eeq
  with the same initial data $u_0$.

  Global and even local solvability of this CP for (\ref{R5})
  is a difficult open problem. Here,  ${\bf A}(u)=|D^2 u|$ is a
  potential operator in $L^2(\re^2)$ with the inner product denoted by $\langle\cdot,\cdot \rangle$.
   The potential is given by
 $$
  \mbox{$
 \Phi(u) = \int\limits_0^1 \langle u, {\bf A}(\rho u) \rangle \, {\mathrm d} \rho= \frac 13 \int u
 |D^2
 u|.
  $}
 $$
  Hence, equation (\ref{R5}) admits
  two integral identities obtained by multiplication by $u$ and
  $u_t$,
   \beq
   \label{R6}
   \begin{matrix}
    \frac 12 \frac {\mathrm d}{{\mathrm d}t} \int u^2= -\e \int(\D
    u)^2 - \int u |D^2 u|, \smallskip\ssk \smallskip \\
    \int\limits_0^T \int (u_t)^2 = E(u(T))-E(u_0), \quad E(u)=-
    \frac \e 2 \int(\D u)^2 - \frac 13 \int u |D^2 u|.
     \end{matrix}
     \eeq
In particular, writing the last identity as
 \beq
   \label{R61} \mbox{$
    \int\limits_0^T \int (u_t)^2 +
    \frac \e 2 \int(\D u)^2 + \frac 13 \int u |D^2 u|=-E(u_0)
   $}
    \eeq
we see that for a uniform bound on $u_\e \ge 0$ it is necessary to
have the following ``dominated concavity property": for all $t \ge
0$,
 \beq
 \label{DC1}
 |D^2 u_\e|>0 \quad \mbox{in  domains,  where $u_\e \ge 0$ is not small}.
  \eeq
  To support this, as a formal illustration,
we  prove an ``opposite nonexistence" result:

\begin{proposition}
 \label{Pr.Non1}
 Assume that, for smooth enough $u_0$,
  \beq
  \label{DC2}
 \mbox{$
 \int u_0 |D^2 u_0|<0.
  $}
  \eeq
  Then $u_\e(\cdot,t)$ for $\e \ll 1$ is not bounded in $L^2$ for large $t>0$.
   \end{proposition}

 Recall that, for the problem (\ref{R5}), the conservation of the mass (i.e., an $L^1$ uniform estimate)
 is available; see Section \ref{SectCons}. Note that (\ref{DC2})
 is not true in the radial case for decreasing $u_0(r) \not \equiv 0$, since by
 (\ref{81}),
 $$
  \mbox{$
 \int u_0 |D^2 u_0|= \int r u_0 \,
 \frac 1 r \, u_0' u_0''=- \frac 12 \int (u_0')^3 >
 0.
 $}
 $$

\smallskip

 \noi{\em Proof.} Estimating from the second identity in
 (\ref{R6})
  $$
  \mbox{$
  \int u|D^2 u| \le - \frac {3\e}2 \int(\D u)^2 + 3\big[\frac \e
  2\int(\D u_0)^2+ \frac 13 \int u_0 | D^2 u_0| \big]
   $}
   $$ and substituting into the first one yields
    \beq
    \label{DC4}
    \mbox{$
    \frac 12  \frac {\mathrm d}{{\mathrm d}t} \int u^2 \ge  \frac \e 2 \int(\D u)^2 - \frac {3 \e}2
    \int(\D u_0)^2 - \int u_0|D^2 u_0|  \ge -\frac  12  \int u_0|D^2 u_0|
    >0
     $}
      \eeq
      for sufficiently small $\e>0$. Hence, under the hypothesis
      (\ref{DC2}),
       \beq
       \label{DC5}
        \mbox{$
       \int u^2(t) \ge \frac 12 \,\big|\int u_0|D^2 u_0|\big| \, t
 \to
       + \infty
       \quad \mbox{as}
       \quad t \to \infty.
        $} \qed
        \eeq

\smallskip

It seems that the divergence (\ref{DC5}) of $\{u_\e\}$ actually
means that the approximated solution $u(x,t)$ is not global and
must blow-up in finite time. This would imply global nonexistence
of solution of (\ref{R4}) if (\ref{DC2}) violates
 the dominant concavity hypothesis (\ref{DC1}) at $t=0$.
 However, we do not know whether  reasonable data satisfying
 (\ref{DC2}) actually exist.
 For instance, the standard profiles in (\ref{zz1}) do not obey
  (\ref{DC2}) (since they correspond to uniformly bounded $L^2$-solutions of (\ref{R4})).

In general,
 identities (\ref{R6}) cannot provide us with
estimates that are sufficient for passing to the limit as $\e \to
0$, so extra difficult analysis is necessary. The main difficulty
is that the Hessian potential $\Phi(u)$ is not definite in the
present functional setting and the operator ${\bf A}(u)$ is not
coercive in the class of not fully concave functions. To avoid
such a difficulty, another uniformly parabolic $\e$-regularization
may be considered useful such as, e.g.,
 \beq
 \label{R10}
 u_\e: \quad  u_t= -  \e \D [(1+u^2) \D u] -|D^2 u|.
  \eeq
Unfortunately, the first operator is not potential in $L^2$, so
deriving integral estimates become more tricky. On the other hand,
using degenerate higher-order $p$-Laplacian operators such as
 \beq
 \label{R11}
 u_\e: \quad  u_t = - \e \D(|\D u|^p \D u) - |D^2 u|
 \eeq
 can be more efficient for $p>1$. Here both operators are potential in
 $L^2(\re^2)$, and moreover, the $p$-Laplacian one is monotone,
 which can simplify derivation of necessary estimates; see Lions' classic
 book \cite[Ch.~2]{LIO}.
 Nevertheless, using various $\e$-regularizations as in
 (\ref{R5}), (\ref{R10}), or (\ref{R11}) does not neglect the
 necessity of the difficult study of boundary layers  as
 $\e \to 0$.

\smallskip

Thus,  according to extended semigroup theory (see
\cite[Ch.~7]{GalGeom}),  the {\em proper} solution of the CP
(\ref{R4}) is given by the limit
 \beq
 \label{R7}
 \mbox{$
 u(x,t) = \lim_{\e \to 0} \, u_\e(x,t).
  $}
 \eeq
 In general, as we have mentioned,
existence (and hence uniqueness) of such limits assumes delicate
studied of $\e$-boundary layers which can occur in the singular
limit $\e \to 0$. In particular, the uniqueness of such a proper
solution would be guaranteed by the fact that the regularized
sequence
 $\{u_\e\}$ does not exhibit $O(1)$ oscillations as $\e \to 0$, so
 (\ref{R7}) does not have different particular limits along
 different subsequences $\{\e_k\} \to 0$.
 Another important aspect is to show that the proper solution does
 not depend on the character of the $\e$-regularizations applied.
 Such a strong uniqueness result is known for the second-order
 parabolic problems \cite[Ch.~6,7]{GalGeom} and is based on the
 MP. For higher-order PDEs, all such uniqueness problems are
 entirely open excluding, possibly, some very special kind of
 equations.
 We hope that the potential properties of the Hessian operator
$|D^2 u|$ and hence identities like (\ref{R6}) can help for
passing to the limit in (\ref{R5}) or other regularized PDEs and
to avoid studying in full generality  difficult singular boundary
layers. These questions remain open.





 \subsection{Riemann's problems: a unique solution via a formal asymptotic series}

 Let us continue to study the passage to the limit $\e \to 0$ in the regularized
problem
   (\ref{R5}). We assume
that the origin 0 belongs to the boundary of ${\rm supp} \, u_0$
and
 \beq
 \label{u1}
 u_0(X)= O (\|X\|^4) \quad \mbox{as} \quad X \to 0 \quad \big(X=(x,y)^T\big).
 \eeq
 This class of data specifies  a kind of a Riemann's problem under consideration (with
 given type of singular transition to 0).
 Then we perform the following scaling in (\ref{R5}):
  \beq
  \label{u2}
  u(X,t)= \e v_\e(\zeta,t), \quad \zeta=X/\e^{\frac 14},
   \eeq
   so that $v_\e(\zeta,t)$ solves an $\e$-independent equation
  with $\e$-dependent data,
 \beq
 \label{u3}
  \mbox{$
 v_\e: \quad  v_t= {\bf A}_0(v)\equiv - \D^2 w - |D^2 w|, \quad v_{0\e}(\zeta) = \frac 1 \e \, u_0(\zeta \e^{\frac 14}).
  $}
   \eeq
According to (\ref{u1}), we assume that there exists a finite limit on any compact subset
 \beq
 \label{u4}
 v_{0\e}(\zeta) \to v_0(\zeta) \quad \mbox{as} \quad \e \to 0,
  \eeq
  where, without loss of generality, by $v_0(\zeta)$ we mean a fourth-degree polynomial.

As usual in asymptotic expansion theory (see e.g. Il'in
\cite{Il92}), the crucial is the first nonlinear step, where we
find the first approximation $V_0(\zeta,t)$ satisfying
  the uniformly parabolic PDE
   \beq
   \label{u5}
   V_0: \quad V_t= {\bf A}_0(V), \quad V(\zeta,0)= v_0(\zeta).
    \eeq
It can be shown by classic parabolic theory \cite{EidSys, Fr} that
the fourth-degree growth of $v_0(\zeta)$ as $\zeta \to \infty$
guarantees at least local existence and  uniqueness of $V_0$.

As the next step, we define the second term of approximation, $v=
V_0 + w_0$, where $w_0$ solves the linearized problem
  \beq
  \label{V13}
w_0: \quad  w_t= {\bf A}_0'(V_0)w,
\quad w(0,y)=v_{0\e}(y)-v_0(y),
\quad \mbox{etc.},
 \eeq
 again  checking that this linear parabolic problem is
 well-posed by classic theory \cite{EidSys, Fr}.

Eventually, this means that we formally express the solution via
the asymptotic series
 \beq
 \label{V14}
  \mbox{$
 v_\e(y,t) = V_0(y,t)+ \sum_{j \ge 0}w_j(y,t;\e),
 $}
  \eeq
  where each term $w_k$ for $k \ge 1$, is obtained by linearization on the
  previous member, by setting
   $$
    \mbox{$
   V_k= V_{k-1}+w_k \equiv V_0+ \sum_{j \le
  k}w_j,
   $}
    $$
     which gives for $w_k$ a non-homogeneous linear parabolic problem
 \beq
  \label{V15}
   \mbox{$
  w_k: \quad w_t= {\bf A}_0'(V_{k-1})w - \big[(V_{k-1})_t -{\bf A}_0(V_{k-1})\big], \quad
w(0,y)=0,
 $}
  \eeq
  with similar assumptions on the well-posedness.

Thus, (\ref{V14}) gives a unique formal representation of the
solution. In asymptotic expansion theory, the convergence of such
series and passing to the limit $\e \to 0$ are often extremely difficult even for
lower-order PDEs, where the rate of convergence or asymptotics are
also hardly understandable. For instance,  the asymptotic expansion for
the classic Burgers' equation
 \beq
 \label{Bur1}
 u_t + u u_x = \e u_{xx}
  \eeq
contains $\ln \e$ terms (close to shock waves), and a technically
hard proof uses multiple reductions of (\ref{Bur1}) to the linear
heat equation, which is illusive for our M-A PDEs; see
\cite{Il92}. For practical reasons, it is important that, as an
intrinsic feature of asymptotic series, each term in (\ref{V14})
(including the first and the simplest one $w_0$) reflects the
actual rate of convergence of $u_\e$ given by (\ref{R7}) as $\e
\to 0$ to the proper solution (\ref{R7}).





\subsection{Conservation of mass}
 \label{SectCons}

In view of (\ref{R7}), this is  straightforward, since
 \beq
 \label{d11}
  |D^2 u| = {\rm div} \, {\bf V}, \quad
  \mbox{where}
  \quad {\bf V} = \mbox{$\frac 12$} \, \left[
   \begin{matrix}
   u_x u_{yy} - u_y u_{xy}\\
  u_y u_{xx}  - u_x u_{xy}
   \end{matrix}
   \right] .
    \eeq
   Therefore, given initial data $u_0 \in L^1(\re)$, integrating (\ref{R5})
  over $\re^2$, we obtain that
 \beq
 \label{ss1}
 \mbox{$
  \int\limits_{\re^2} u_\e(x,t) \, {\mathrm d} x = \int\limits_{\re^2} u_0(x) \, {\mathrm d}
  x \quad \mbox{for} \quad t >0,
   $}
    \eeq
    meaning the {\em conservation of mass} of this regularized  M-A flow.
In this sense, (\ref{R5}) looks like a fully nonlinear version of
the quadratic PME from filtration theory
 $$
 u_t = \D u^2 \equiv  {\rm div} \, {\rm grad} \, u^2,
  $$
  but corresponds to more nonlinear and complicated mathematics
  without standard MP and accompanying
  barrier techniques.


 \section{Examples of blow-up for a fourth-order M-A equation with $-|D^4 u|$}
  \label{Sect4ord}

\subsection{On  derivation of the higher-order radial M-A model}

Our intention is to show that such models make sense and exhibit
correct local parabolic well-posedness. Namely, we introduce
radial models related to the fourth-order M-A equation
 \beq
 \label{ma4}
 u_t = - |D^4 u| + u^3
 \quad \mbox{in} \quad \re^2\times
 \re,
  \eeq
  with the catalecticant determinant $|D^4 u|$  given by
 \beq
 \label{F41}
  {\rm det} D^4 u \equiv {\rm det} \left[ \begin{matrix} u_{xxxx}\,\,
 u_{xxxy}\,\,
 u_{xxyy} \cr
  u_{xxxy}\,\, u_{xxyy}\,\, u_{xyyy} \cr
  u_{xxyy}\,\, u_{xyyy}\,\, u_{yyyy} \end{matrix} \right],
  \eeq
  which plays an important role in the theory of quartic forms.
  For instance, each such form in two variables can be expressed {via} a
  sum of three fourth powers of linear forms and {via} two powers,
  provided that ${\rm det} \, D^4 u=0$; see \cite{Har92}.

 It then follows from (\ref{F41}) that
  \beq
  \label{D4}
   \begin{matrix}
  |D^4 u|= u_{xxxx}u_{yyyy}u_{xxyy}+2 u_{xxxy}u_{xyyy}u_{xxyy} \ssk\ssk\\
 - \,(u_{xxyy})^2-  (u_{xyyy})^2 u_{xxxx}-   (u_{xxxy})^2
 u_{yyyy}.
\end{matrix}
 \eeq
 In particular,
 for radial functions $u=u(r)$, we have
  \beq
  \label{DD1}
   \begin{matrix}
  u_{xxxx} = \frac {x^4}{r^4} \, u^{(4)} + \frac{6x^2y^2}{r^5} \,
  u''' + \frac{3(y^4-4x^2y^2)}{r^6} \, u''
  + \frac{3(-y^4+4x^2y^2)}{r^7} \, u', \ssk \\
u_{yyyy} = \frac {y^4}{r^4} \, u^{(4)} + \frac{6x^2y^2}{r^5} \,
  u''' + \frac{3(x^4-4x^2y^2)}{r^6} \, u''
  + \frac{3(-x^4+4x^2y^2)}{r^7} \, u', \ssk \\
u_{xxyy} = \frac {x^2 y^2}{r^4} \, u^{(4)} +
\frac{x^4+y^4-4x^2y^2}{r^5} \,
  u''' + \frac{-2(x^4+y^4)+ 11 x^2y^2}{r^6} \, u''
  + \frac{2(x^4+y^4)- 11 x^2y^2}{r^7} \, u', \ssk \\
u_{xxxy} = \frac {x^3 y}{r^4} \, u^{(4)} + \frac{3(x y^3- x^3 y)
}{r^5} \,
  u''' + \frac{3(2 x^3 y - 3  x y^3)}{r^6} \, u''
  + \frac{3(-2 x^3 y+ 3 x y^3)}{r^7} \, u', \ssk \\
u_{xyyy} = \frac {x y^3}{r^4} \, u^{(4)} + \frac{3(x^3 y- x y
^3)}{r^5} \,
  u''' + \frac{3(2 x y^3 - 3  x^3 y)}{r^6} \, u''
  + \frac{3(-2 x y^3+ 3 x^3 y)}{r^7} \, u'.
   \end{matrix}
 \eeq
Overall, this leads to the following  radial model that is the
simplest in such flows:


\subsection{Simplest well-posed radial fourth-order M-A type model with regional blow-up}

Using (\ref{DD1}), by balancing and mutual cancellation of the
most singular terms as $r \to 0$, we introduce the following model
radial parabolic equation
  associated with the M-A flow (\ref{ma4}):
  \beq
  \label{ma50}
   \mbox{$
  u_t= - \frac 1{r^2} \, (u_{rrr})^2 u_{rrrr} + u^3
  \quad \mbox{in} \quad \re_+ \times \re_+.
   $}
 \eeq
 Looking at it as a parabolic equation, we pose at the origin the
 symmetry (regularity) conditions
  \beq
  \label{sr10}
  u_r(0,t)=u_{rrr}(0,t)=0.
 \eeq
Obviously, (\ref{ma50}) is uniformly parabolic with smooth
(analytic) solutions in any domain of non-degeneracy $\{u_{rrr}
\not = 0\}$. In particular, checking the regularity of the
operator in (\ref{ma50}) and passing to the limit $r \to 0$ yields
 \beq
 \label{rr10}
  \mbox{$
  - \frac 1{r^2}\, (u_{rrr})^2 u_{rrrr} \to -(u_{rrrr})^3,
   $}
   \eeq
 so that, at the origin, the differential operator is
 non-degenerate and regular if $u_{rrrr} \not = 0$.
 Thus, regardless the  degeneracy of the equation (\ref{ma50}),
 this radial version of fourth-order M-A flows is well-posed, at
 least, locally in time.

Describing  blow-up patterns and  looking for blow-up similarity
solutions of (\ref{ma50}),
 \beq
 \label{rr30}
  \mbox{$
 u_{\rm S}(r,t)= \frac 1{\sqrt{T-t}} \,
    f(r),
     $}
  \eeq
 we obtain the ODE
  \beq
  \label{rr40}
   \mbox{$
   {\bf A}(f) \equiv
  - \frac 1 {r^2} \, (f''')^2 f^{(4)} + f^3 - \frac 12 \, f = 0
   \quad \mbox{for} \quad r>0,
   $}
   \eeq
   with the  symmetry conditions generated by (\ref{sr10}) at the
   origin,
 \beq
    \label{symm1L}
    f'(0)= f'''(0)=0.
     \eeq

\begin{figure}
\centering
\includegraphics[scale=0.65]{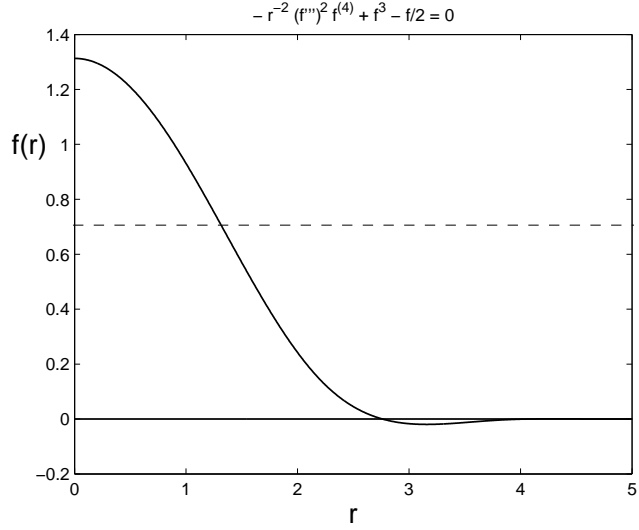}
\caption{\small A unique compactly supported nonnegative  solution
of the ODE (\ref{rr40}).}
 \label{F4SS}
\end{figure}

   Indeed, (\ref{rr40}) is a difficult fourth-order ODE with non-monotone, non-autonomous, and
   non-potential operators.
Any (at least) 3D phase-plane analysis or shooting arguments are
rather difficult. Instead, and this is much easier nowadays, we
are going to
  check existence (and uniqueness) of
the solutions using reliable numerical methods of an enhanced
accuracy. For instance, Figure \ref{F4SS} demonstrates the
corresponding compactly supported solution of (\ref{rr40}). We did
not detect other P-type profiles that have more than one
oscillation about the constant equilibrium
 $$
  \mbox{$
 f_0= \frac 1{\sqrt 2}.
  $}
  $$

As usual, stability of the similarity blow-up (\ref{rr30}) is
studied in terms of the rescaled solution
 \beq
 \label{w1}
 w(x,\t) = \sqrt{T-t}\,\, u(r,t), \quad \t= - \ln(T-t) \to + \infty,
  \eeq
 where $w$ solves the rescaled parabolic equation with the same
elliptic operator
 \beq
 \label{w2}
  \mbox{$
 w_\t ={\bf A}(w) \equiv - \frac 1{r^2} \, (w_{rrr})^2 w_{rrrr} + w^3 - \frac
 12 \, w \quad \mbox{in} \quad \re_+ \times \re_+.
  $}
   \eeq
Since ${\bf A}$ is not potential in any suitable metric, so
(\ref{w2}) is not a gradient system, passage to the limit as $\t
\to +\infty$ to show stabilization to the stationary profile
$f(r)$ represents a difficult open problem. Note that the main
operator  $ - \frac 1{r^2} \, (w_{rrr})^2 w_{rrrr}$
 is gradient and
 admits
multiplication by $r^2 \, w_{rr\t}$ in $L^2(\re)$,
 $$
  \mbox{$
 -\int (w''')^2 w'''' w_\t '' = \frac 1{12}\,
 \frac{\mathrm d}{{\mathrm d}\t} \, \int
 (w''')^4 \quad ('=D_r),
  $}
  $$
 but $w^3$ (and also $w$) does
not.

\subsection{$p=3$: on oscillatory structure close to interfaces}

\begin{figure}
\centering
\includegraphics[scale=0.65]{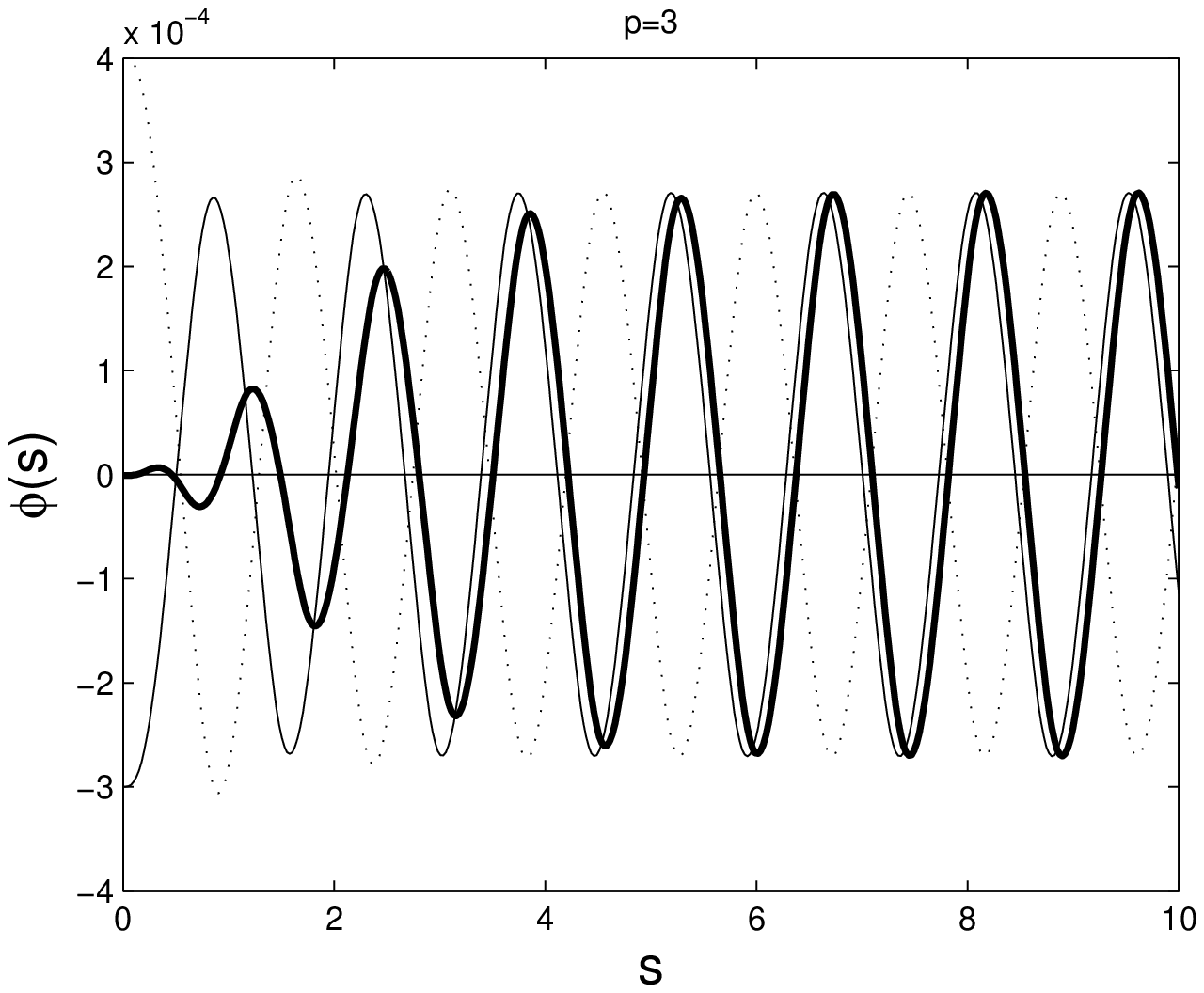}
\caption{\small Stabilization to a stable periodic orbit
$\varphi_*(s)$ for the ODE (\ref{o20}), $\l_0=1$.}
 \label{Fosc1100}
\end{figure}

In addition, Figure \ref{F4SS} shows that, locally, close to the
finite interface point $r=r_0$, sufficiently smooth solutions of
(\ref{rr40}) are oscillatory. This kind of non-standard behaviour
of solutions of the Cauchy problem deserves a more detailed
analysis.
 To describe this, we
introduce an extra scaling by setting
 \beq
 \label{o10}
  \mbox{$
 f(r) = (r_0-r)^\g \varphi(s), \quad s= \ln(r_0-r), \quad
  \mbox{where} \quad  \g =3.
  $}
  \eeq
Substituting (\ref{o10}) into (\ref{rr40}) and neglecting the
higher-degree term $f^3$ for $r \approx r_0^-$, we obtain the
following equation for the {\em oscillatory component}
 $\varphi(s)$:
  \beq
  \label{o20}
   \mbox{$
   (f''')^2 f^{(4)}=- \l_0 f \quad \Longrightarrow
   \quad
  (P_3(\varphi))^2 P_4(\varphi)= -\l_0
  \varphi, \quad \l_0= \frac 12 \, r_0^2,
   $}
   \eeq
   where $P_k$ are linear differential polynomials
   obtained by the recursion procedure (see
   \cite[p.~140]{GSVR})
 $$
 \begin{matrix}
P_1(\varphi)=\vp'+ \g \vp, \quad P_2(\vp)= \vp''+(2\g-1) \vp'+
\g(\g-1)\vp, \smallskip\smallskip \\
  P_3(\varphi)=\vp'''+ 3(\g-1) \vp'' + (3 \g^2- 6 \g +2) \vp'
 + \g(\g-1)(\g-2) \vp, \smallskip\smallskip \\
  P_4(\varphi)=\vp^{(4)}+ 2(2\g-3) \vp''' + (6 \g^2- 18 \g +11)
  \vp'' \smallskip\smallskip \\
   +\, 2(2 \g^3 -9 \g^2 + 11 \g -3) \vp'
 + \g(\g-1)(\g-2)(\g-3) \vp.
 \end{matrix}
 $$
In Figure \ref{Fosc1100}, we show the typical behaviour of
solutions of the second ODE in  (\ref{o20}) demonstrating a fast
stabilization to the unique stable periodic solution. According to
(\ref{o10}), this periodic orbit $\varphi_*(s)$
describes the generic character of oscillations of solutions of
the ODE (\ref{rr40}).
 Periodic solutions for {\em semilinear}
 ODEs of the type (\ref{o20}) are already known for the third-
 \cite[\S~7]{Gl4} and fifth-order operators (with  $P_5(\varphi)$)
 \cite[\S~6,\,12]{GBl6}.
 The quasilinear equation (\ref{o20}) is more difficult,
  and
 existence and uniqueness of $\varphi_*(s)$ remain open.




As a key application of the expansion (\ref{o10}),
 we note that this shows that the ODE
 (\ref{rr40}) generates a 2D asymptotic bundle of solutions
  close to interfaces
  \beq
  \label{2d1}
  f(r)= (r_0-r)^3 \varphi(s+s_0)+... \quad \mbox{as \,\, $r \to r_0^-$},
  \eeq
  with two parameters, $r_0>0$ and $s_0 \in \re$ as the phase
  shift in the periodic orbit. The 2D bundle perfectly suits
  shooting also {\em two} boundary conditions (\ref{symm1L}),
  though a proper topology of shooting needs and deserves  extra analysis.

\subsection{Structures of single point blow-up for $p>3$}

As usual, single point blow-up occurs for $p>3$ in the M-A PDE
with other power nonlinearities,
  \beq
  \label{ma50L}
   \mbox{$
  u_t= - \frac 1{r^2} \, (u_{rrr})^2 u_{rrrr} + |u|^{p-1}u
  \quad \mbox{in} \quad \re_+ \times \re_+,
   $}
 \eeq
 where we again pose the symmetry conditions (\ref{sr10}).
The blow-up similarity solutions are now
 \beq
 \label{rr30L}
  \mbox{$
 u_{\rm S}(r,t)= (T-t)^{-\frac 1{p-1}} f(z), \quad z= r/(T-t)^\b,
 \quad \mbox{where}
 \quad \b= \frac {p-3}{12(p-1)},
     $}
  \eeq
 and $f$ solves the ODE
  \beq
  \label{rr40L}
   \mbox{$
   {\bf A}(f) \equiv
  - \frac 1 {z^2} \, (f''')^2 f^{(4)}- \b f' z - \frac 1{p-1} \, f+|f|^{p-1} f = 0
   \quad \mbox{for} \quad z>0,
   $}
   \eeq
   with the same symmetry conditions at $z=0$, (\ref{symm1L}).
   For $p=3$, (\ref{rr40L}) yields the simpler autonomous equation
   (\ref{rr40}) for regional blow-up.
For $p>3$, this ODE is more difficult, and, following the results
in Section \ref{S2}, we expect to have blow-up profiles of P, Q
and, possibly, S-type.

\begin{figure}
\centering
\includegraphics[scale=0.7]{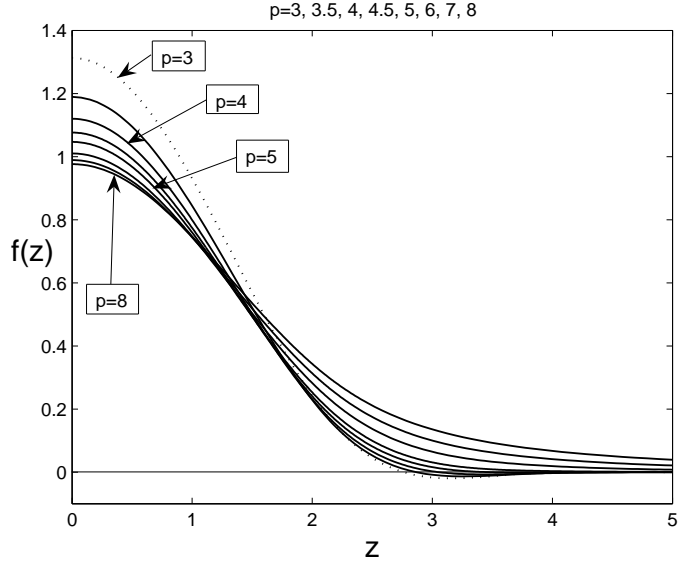}
\caption{\small Single point blow-up profile satisfying  the ODE
(\ref{rr40L}) for $p \in(3,8]$.}
 \label{F4LS}
\end{figure}

    Figure \ref{F4LS}
demonstrates the first profiles of P-type for various $p \in
[3,8]$. The dotted line corresponds to the regional blow-up
profiles, $p=3$, from Figure \ref{F4SS}. In particular, this
clearly shows the continuity of the (``homotopic") deformation of
solutions of (\ref{rr40L}) as $p \to 3^+$.

\begin{figure}
\centering
\subfigure[enlarged, $10^{-3}$]{
\includegraphics[scale=0.52]{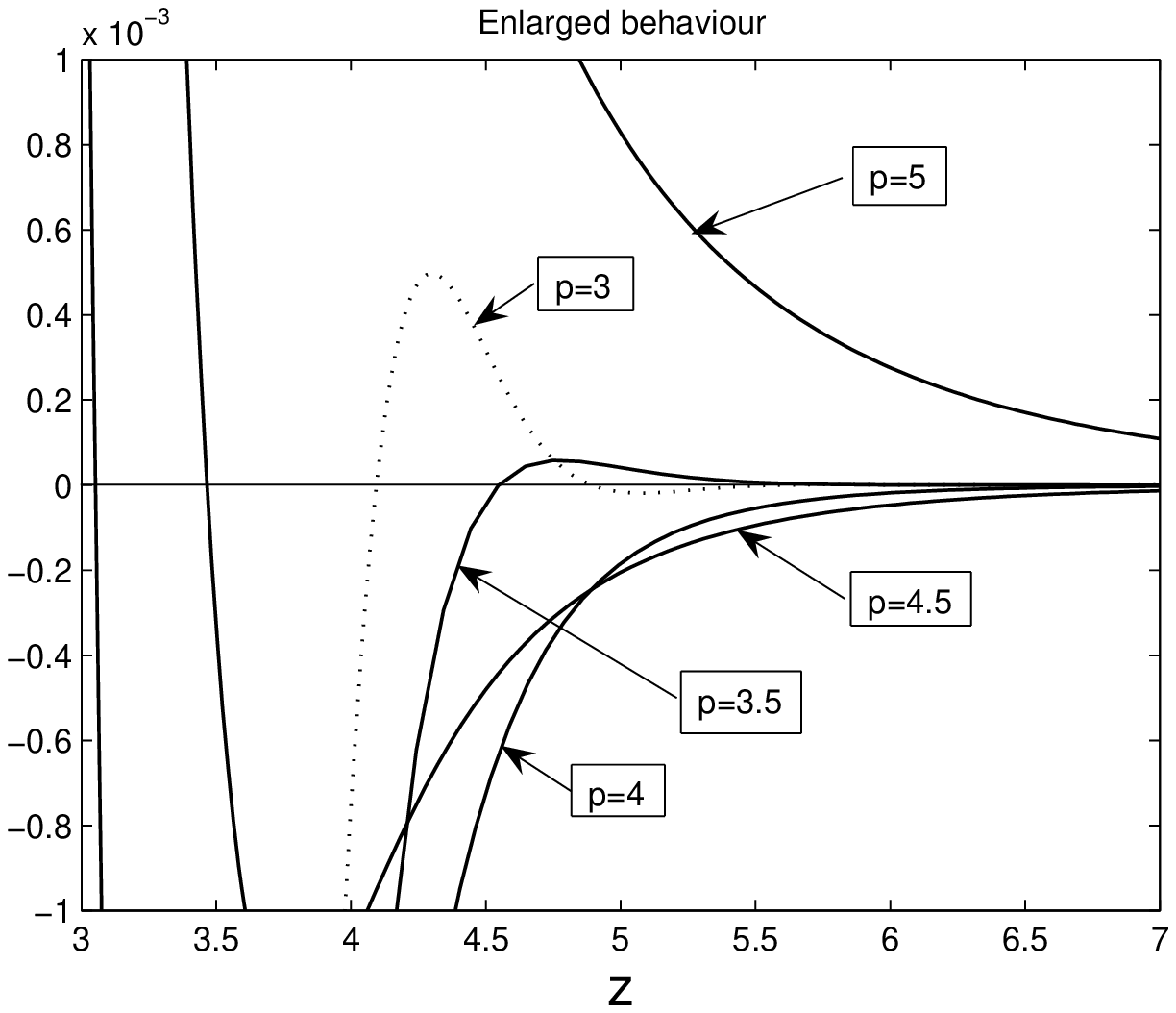} 
}
 \subfigure[more enlarged, $10^{-4}$]{
\includegraphics[scale=0.52]{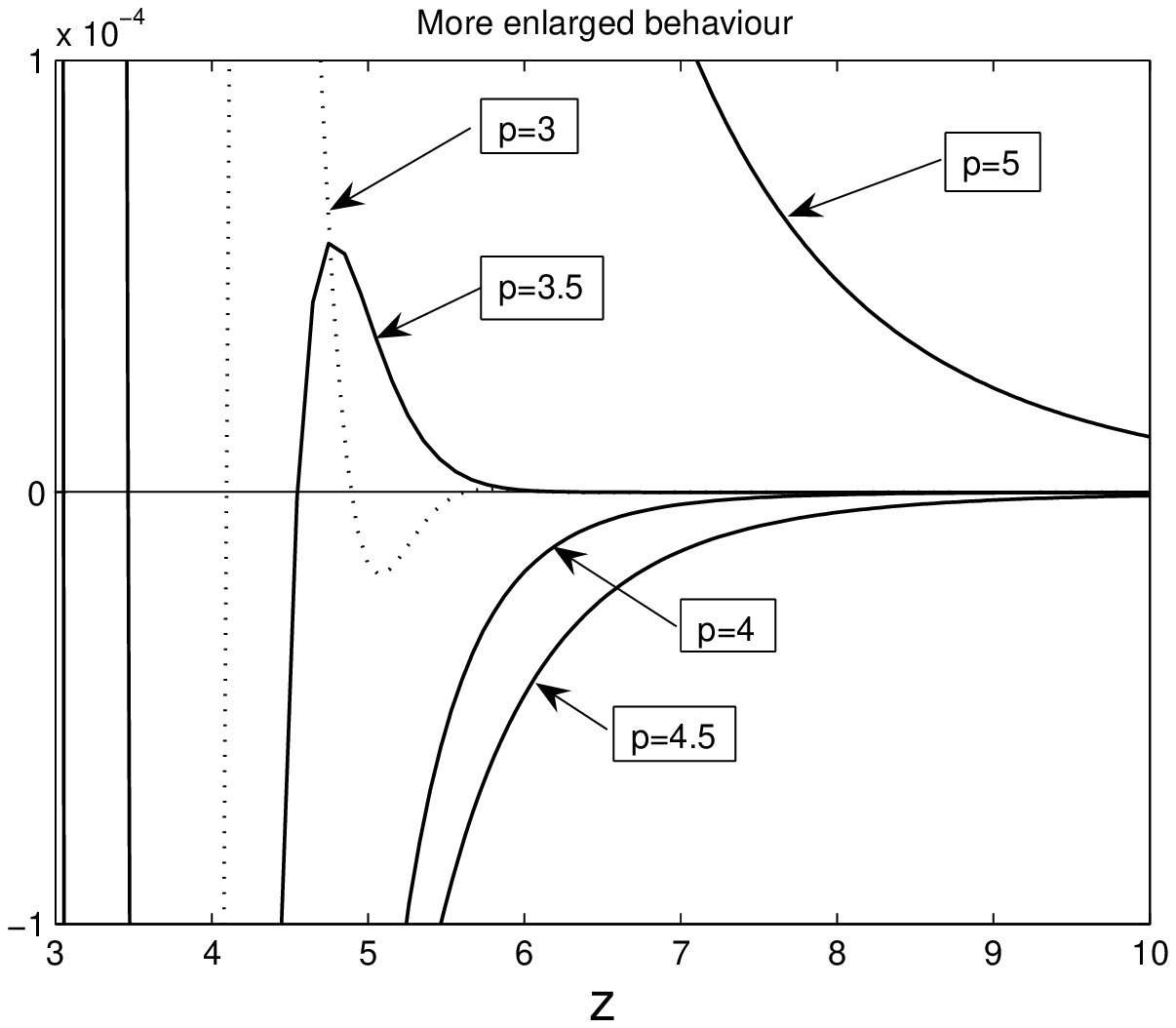} 
}
   \vskip -.3cm
    \caption{\small Non-oscillatory behaviour of profiles
    from Figure \ref{F4LS} for $p>3$.}
   \vskip -.3cm
 \label{FFss1}
\end{figure}

It is key to observe that the profiles for $p>3$ have infinite
interface and, moreover,
 \beq
 \label{mm1}
 f(z) > 0 \,\,\, \mbox{in} \,\,\, \re_+ \,\,\, \mbox{for $p$
 larger than, about, 5.}
 \eeq
 For smaller $p > 3$, the profiles continue to change sign as
 for $p=3$, as the continuity in $p$ suggests.
 In Figure \ref{FFss1}, we show the enlarged behaviour
 of the profiles from Figure \ref{F4LS} in the domains, where these
 are sufficiently small. In both Figures (a) and (b), the profile
 for $p=3.5$ has two zeros only.

It  follows from  (\ref{rr40L}) that the non-oscillatory profiles
$f(z)$ have the asymptotic behaviour governed by two linear terms
(other nonlinear ones are negligible on such asymptotics),
 \beq
 \label{as1}
  \mbox{$
- \b f' z - \frac 1{p-1} \, f=0
 \quad \Longrightarrow \quad f(z) = C z^{\nu}+... \, , \quad
 \nu=-\frac{12}{p-3}<0,
  $}
  \eeq
  where $C=C(p) \not= 0$ for a.e. $p>3$.
 The full 2D bundle of such non-oscillatory asymptotics includes an
 essentially ``non-analytic"  term of a typical centre manifold nature
 \cite{Perko}:
  \beq
  \label{nn1}
   \mbox{$
  f(z) = C z^{\nu}+...+ C_1 {\mathrm e}^{-b_0 z^\mu}+... \, ,
  \quad \mbox{where} \quad  \mu= \frac {11 p-9}{p-3}>0,
   $}
   \eeq
  $C_1 \in \re$ is the second parameter,
 and
    $b_0=b_0(p)>0$ is a constant that can be easily computed.
    The expansion (\ref{nn1}) can be justified by rather technical
    application of standard fixed point theorems in a weighted
    spaces of continuous functions defined for large $z \gg 1$.

The full 2D bundle (\ref{nn1}) poses a well-balanced matching
problem to satisfy {\em two} symmetry conditions at the origin
(\ref{symm1L}) for any $p>3$, provided that $C \not = 0$.

On the other hand,
   Figure \ref{FFss1}(b)
  clearly shows that $C(p)$ changes sign  at some
   \beq
   \label{pp1}
   \hat p_1 \in (4.5,5), \quad \mbox{and} \quad C(\hat p_1)=0.
    \eeq
Then, at $p=p_1$, using the oscillatory analysis presented below,
we expect that the corresponding $f(z)$ is compactly supported.

Moreover,  in view of the oscillatory behaviour near interfaces
(cf. (\ref{o10}) for $p=3$), we expect that there exists a
monotone decreasing sequence of such critical exponents
 \beq
 \label{ss11}
 \{\hat p_k\}_{k \ge 1} \to 3^+ \,\,\, \mbox{as} \,\,\,k \to
 \iy,
 \quad \mbox{such that} \quad C(\hat
 p_k)=0,
  \eeq
so that $\hat p_1$ in (\ref{pp1}) is just the first, maximal one.
Possibly, such a mixture of compactly supported and non-compactly
supported profiles $f(z)$ for $p>3$ gets more complicated (see
further comments below on oscillatory character of
finite-interface solutions).

In the cases (\ref{ss11}), the single point blow-up profile can be
compactly supported, so we need to describe its local behaviour
near the interface, which turns out to be different from that for
$p=3$ studied above.

\subsection{Finite interfaces: on  oscillatory structures for $p>3$}

 The oscillatory structure of solutions is given by two first
terms in (\ref{rr40L}), and hence have a different form, than in
(\ref{o10}),
 \beq
 \label{o10L}
  \mbox{$
 f(z) = (z_0-z)^\g \varphi(s), \quad s= \ln(z_0-z), \quad
  \mbox{where} \quad  \g =\frac 92.
  $}
  \eeq
Substituting into (\ref{rr40L}) and neglecting other higher-degree
terms, yields for
 $\varphi(s)$ the ODE
  \beq
  \label{o20L}
   \mbox{$
  (P_3(\varphi))^2 P_4(\varphi)= -
  P_1(\varphi)\equiv -( \varphi' + \g \vp),
   $}
   \eeq
   where we have scaled out the constant multiplier
     $\b z_0^3>0
     $ on the right-hand side. This ODE is even more difficult
     than (\ref{o20}), so we again rely on careful numerics.

\begin{figure}
\centering
\includegraphics[scale=0.65]{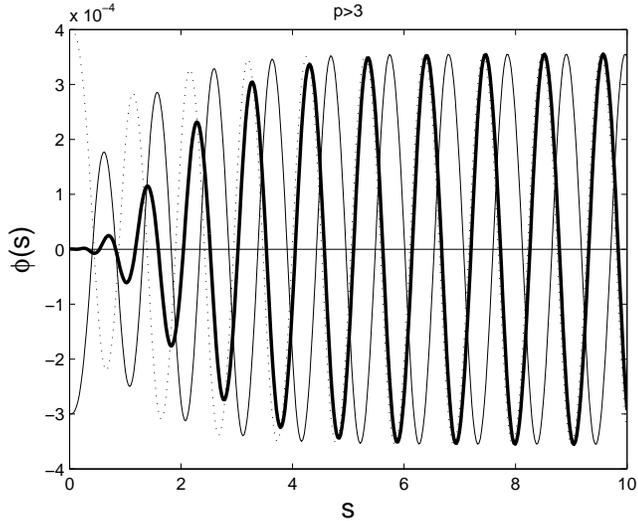}
\caption{\small Stabilization to a stable periodic orbit
$\varphi_*(s)$ for the ODE (\ref{o20L}).}
 \label{Fosc1100L}
\end{figure}

In Figure \ref{Fosc1100L}, we show the typical behaviour of
solutions of (\ref{o20L}) demonstrating a fast stabilization to a
unique and stable periodic solution. According to (\ref{o10L}),
this periodic orbit $\varphi_*(s)$  describes the generic
character of oscillations of solutions of the ODE (\ref{rr40L}) in
the critical case (\ref{pp1}).

Since here the periodic orbit is stable as $s \to +\infty$, it is
unstable as $s \to -\infty$ (in the direction towards the
interface at $s=-\infty$ according to (\ref{o10L})), and moreover
the stable manifold of $\varphi_*(s)$  as $s \to -\infty$ consists
of the solution itself up to shifting. Therefore,  the full
equation (\ref{rr40L}) admits precisely
 \beq
 \label{2dL}
 2D \,\,\,\mbox{bundle of small solutions, with parameters $z_0>0$
 and $s_0 \in \re$},
  \eeq
  where $s_0$ is again the translation in the oscillatory periodic
  component $\vp(s+s_0)$. Therefore, shooting via 2D bundle
  precisely {\em two} symmetry conditions at the origin
  (\ref{symm1L}) represents a well-posed problem, which can admit
  solutions, and possibly a countable set of these.
  Nevertheless, Figures \ref{FFss1}(a) and (b) justify that the
  actual behaviour is governed by the non-compactly supported 2D
  bundle (\ref{nn1}) for a.e. $p>3$, and we do not know any
  mathematical reason why the finite-interface bundle (also 2D)
 fails to be applied for $p>3$ a.e. and not only for $p=p_k$.

\subsection{On multiplicity of solutions by branching}

\begin{figure}
\centering
\includegraphics[scale=0.65]{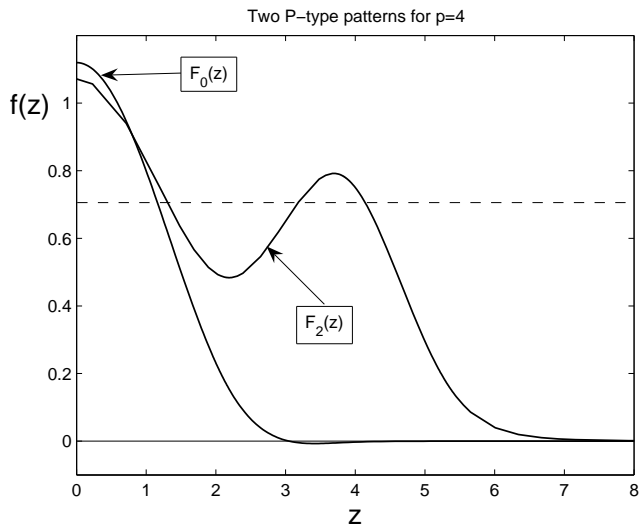}
\caption{\small Two P-type profiles of the ODE (\ref{rr40L}) for
$p=4$.}
 \label{F4LSNN}
\end{figure}

 Figure \ref{F4LSNN} shows two
 P-type profiles, $F_1(z)$ and  $F_2(z)$. This poses a difficult open problem
 on multiplicity of solutions for $p>3$ (and also for $p=3$).
Recall that operators in both  equations (\ref{rr40L}) and
(\ref{rr40}) are not potential, so we cannot rely on well
developed theory of multiplicity for variational problems; see
e.g., \cite[\S~57]{KrasZ}.

Nevertheless, we will rely on variational theory by introducing a
family of approximating operators
 \beq
 \label{aa1}
  \mbox{$
 {\bf A}_\mu(f)= \big[- \frac \mu{z^2}\, (f''')^2 +
 \mu-1\big]f^{(4)} - \mu \b f'z- \frac 1{p-1}\, f+
 |f|^{p-1}f, \quad \mu \in[0,1].
  $}
  \eeq
Then ${\bf A}_1={\bf A}$, while for $\mu=0$,
 \beq
 \label{aa2}
  \mbox{$
 {\bf A}_0(f)= - f^{(4)} -\frac 1{p-1}\, f+ |f|^{p-1}f,
  $}
  \eeq
  we obtain a standard non-coercive variational operator. The
  corresponding functional for $f \in H^2(-L,L)$ on a fixed interval with $L \gg
  1$,
   \beq
   \label{aa3}
   \mbox{$
   \Phi_0(f)= - \frac 12 \int(f'')^2 - \frac 1{2(p-1)}\, \int f^2  + \frac 1{p+1}\, \int
   |f|^{p+1}
 $}
  \eeq
has at least a countable set of different critical points
$\{f_k^{(0)}, \, k=0,1,2,...\}$; see \cite{PohFM}.

Thus we arrive at the branching problem from profiles $f_k^{(0)}$
at the branching point $\mu=0$, which leads to classic branching
theory; see \cite[Ch.~6]{KrasZ} and \cite{VainbergTr}. In general
in the present ODE setting, where the linearized operator
 for $\mu=0$,
  \beq
  \label{aa4}
  {\bf A}_0'(f)Y= - Y^{(4)}+ p|f|^{p-1}Y,
  \eeq
has a 1D kernel, branching theory \cite[p.~401]{Deim} suggests
that each member $f_k^{(0)}$ generates a continuous branch
$\{f_k^{(\mu)}\}$ at $\mu=0$. The global continuation of those
branches up to $\mu=1$ remains a difficult open problem, that in
the present ODE case admits an effective numerical treatment.



\begin{small}
\begin{appendix}

\section*{Appendix A: On various Monge--Amp\`ere models and
application}

 \label{ApA}
 \setcounter{section}{1}
\setcounter{equation}{0}

\subsection*{Second-order Monge--Amp\`ere equations:
 critical exponents and  singularities}

Thus, for a given function $u \in C^2(\ren)$,  $D^2u$ denotes the
corresponding $N \times N$ {\em Hessian} matrix\index{Hessian}
$D^2 u=\|u_{x_ix_j}\|$. As we have mentioned in Introduction,
general parabolic {\em Monge--Amp\`ere} (M-A) {\em equations}
(\ref{MA.1})
 play a key role in many  geometric problems and applications.
Thus,  M-A equations first appeared in  Monge's paper
\cite{Monge81} in 1781,\index{Monge, G.} related to
civil-engineering problem of moving a mass of earth from one
configuration to another in the most economical way. This  was
further studied by
 Appel \cite{Appel87} in 1887
   and
Kantorovich  in the 1940s\cite{Kant42, Kant48};
  see
history in \cite{Feyel03}. Other key problems and M-A applications
are:

 (i) logarithmic 
Gauss and Hessian  curvature flows,

 (ii) the Minkowski problem (1897),

 (iii) the Weyl problem (with  Calabi's related conjecture in complex geometry),
 etc.

 \noi We refer to
 basic mathematical results  in monographs by
 Taylor   \cite[Ch.~14,15]{Tay},  Gilbarg--Trudinger \cite[Ch.~17]{GilTr},
  and Guti\'errez \cite{Gut01}.
For increasing functions $g(s)$,  equation (\ref{MA.1}) is
parabolic if $D^2u(\cdot,t)$ remains positively definite for
$t>0$, assuming that $D^2u_0>0$ for initial data $u_0$. For a
class of lower-order  operators $h(\cdot)$ satisfying necessary
growth estimates, typical monotone increasing, concave
nonlinearities $g$ in the principal operator are
 $$
 \mbox{$
 g(s) = \ln s, \,\,\, g(s) = - \frac 1s, \,\,\,\mbox{and} \,\,\,
g(s) = s^{\frac 1N} \quad \mbox{for} \,\,\, s> 0, $}
 $$
 for which the global-in-time solvability is known.
 For other  $g(s)$ functions with a faster growth as $s \to \infty$,
 local-in-time solutions existing by  standard parabolic
 theory (see {e.g.,} the classic book \cite[p.~320]{Lad68})
  may blow-up in finite time. This  gives  special asymptotic patterns,
  which can
  also be of interest in some geometric applications. For
  instance, the affine normal flow for an initial convex, properly
  embedded, and noncompact hypersurface ${\mathbb L} \subset
  \re^{N+1}$ can be describes by the parabolic equation
   \beq
    \label{Loft1}
  u_t= - ({\rm det}\, D^2)^{-\frac 1{N+2}};
   \eeq
  see \cite{Loft08} as a guide for the history and recent results including finite-time extinction.

\ssk

 The parabolic M-A equation
 $$
 -u_t \, {\det} D^2u = f \quad \mbox{in} \,\,\, Q_T= \O \times (0,T),
 $$
 where $\O \subset \ren$ is a bounded smooth domain, was first introduced in \cite{Kr76};
 see recent related
  references in \cite{GutH01} and more general models like that in \cite{Sch03}.~Note
  the pioneering paper by Hamilton
  \cite{Ham82}\index{Hamilton, R.S.} on the evolution of a metric
  in direction of its Ricci curvature.
 Conditions of the global unique solvability of the   M-A equation
  $$
u_t = ({\rm det} D^2u)^{\frac 1N} + g \quad \mbox{in} \,\,\, Q_T,
 $$
were obtained by Ivochkina and Ladyzhenskaya \cite{Iv94}.
This model corresponds to special curvature flows.

 The logarithmic Gauss curvature flow in terms of  the support
 function (see \cite{Chou00}) is described by the M-A-type
 equation
  $$
  \mbox{$
  u_t = \sqrt{1+|x|^2}\, \ln \bigl({\rm det} D^2 u \bigr) + h
  \quad \mbox{in} \,\,\, \ren \times \re_+,
   $}
  $$
  where $h$ is given. Depending on the initial uniformly convex hypersurfaces, this
   PDE is known to
  admit either a local solution, corresponding to shrinking to a
  point in finite time or a global solution that describes the uniform
  convergence
  to an expanding sphere.

  The relations between K\"ahler--Ricci flows on $M$ ($(M^n,g_0)$
  is a complete K\"ahler manifold)
   $$
   \mbox{$
   \frac{\partial g_0}{\partial t}=-R_{i \bar j}, \quad g_{i \bar
   j}(x,0)=(g_0)_{i \bar j}
   $}
   $$
    and the corresponding parabolic M-A equations
     $$
     \mbox{$
      \frac{\partial u}{\partial t}= \ln \frac{{\rm det}\,((g_0)_{k \bar j}+
      u_{k \bar j})}{{\rm det}\, (g_0)_{k \bar j})}- f_0, \quad
      u(x,0)=0,
      $}
      $$
      is explained in \cite{Chau09}, where further references on
      global solvability, asymptotics,  convergence, and application of these problems can be
      found.

Typical second-order  Hessian operators are known to be potential
and  the corresponding smooth parabolic flows  are gradient
systems. For M-A PDEs, these ideas go back to Bernstein (1910)
\cite{Bern10};
 see also Reilly,
 and \cite{PT68} for a particular
case.
 Potential properties of more general Hessian
operators are described  in \cite{Tso90}.
 Parabolic M-A equations as gradient
flows, and the related questions of the asymptotic behavior of
solutions, were studied, {e.g.,} in \cite{Chou01, Sch03}, where
 further references concerning various types of Gaussian flows can
be found.
 The classical Gauss
 curvature flow describes the deformation of a convex compact
 surface $\Sigma: z=u(x,y,t)$ in $\re^3$ by its  Gauss curvature
and is governed by the PDE
 \beq
 \label{GG.11}
 \mbox{$
 u_t = \frac{{\rm det} D^2u}{(1+|\n u|^2)^{3/2}},
 $}
 \eeq
 which is uniformly parabolic  on strictly convex
 solutions.
Singularity formation phenomena for (\ref{GG.11}) appear if the
initial surface $\Sigma$ has flat sides, where the curvature
becomes zero and the equation  degenerate. This leads to an FBP
for (\ref{GG.11}) with the unknown domain of singularity,
$\{(x,y): \, u(x,y,t)=0\}$, and specific regularity properties;
see \cite{DH99, DL03} and references therein. Alternatively,
finite-time formation of {\em non-smooth} free boundaries with
flat parts is a typical phenomenon for blow-up solutions of the
reaction-diffusion PDEs (see the first equation in (\ref{RD.991})
below) {via} extended semigroup theory. In 1D, optimal regularity
of such $C^{1,1}$-interfaces is well understood
\cite[Ch.~5]{GalGeom}. For $N>1$, the regularity problem remains
essentially
 {open};
   see some estimates and examples in
\cite[p.~151]{GalGeom}. As a formal extension, note that a
nontrivial ($u(x,t) \not \equiv 0$) proper convex solution exists
for such PDEs with an arbitrarily strong (as $u \to 0$) absorption
term, {e.g.,}
 $$ 
 \mbox{$
 u_t = \frac{{\rm det} D^2u}{(1+|\n u|^2)^{ 3/2}} - {\mathrm
 e}^{1/u},
 $}
 $$ 
 where a similar FBP  occurs. Therefore, the parabolic operator of
 the Gauss curvature flow is extremely powerful, in the sense that it prevents a {\em complete
extinction} ({i.e.,} $u(x, t) \equiv 0$ for arbitrarily small
$t>0$; this can happen for many other parabolic PDEs). For any
initial data with flat sides, $\{(x,y): \,u_0(x,y)=0\} \not =
\emptyset$, this proper solution of the FBP can be constructed by
regular approximations of the equations and initial data by
replacing
 $
  \mbox{$
 {\mathrm e}^{1/u} \mapsto {\rm min}\{\frac u \e, \,
{\mathrm e}^{ 1/u}\}, \quad \mbox{with} \quad \e
> 0,
 $}
 $
  (a
uniformly Lipschitz continuous approximation), and  $u_0 \mapsto
u_0+ \e$. Uniform {\em a priori} estimates for $\{u_\e\}$ are
obtained by local (near the interface) comparison with 1D TW
solutions or other radial sub- and super-solutions,
\cite[Ch.~7]{GalGeom}.

In connection with other PDEs, let us also mention an unusual
hyperbolic M-A  equation
 $$
 u u_{tt} + u_{tt}u_{xx}-(u_{xt})^2+1=0,
  $$
  which occurs in 1D hyperbolic mean curvature flows for closed
  curves $F_{tt}=k {\bf N}- \n p$ on $S^1 \times [0,T)$, $\n p=
  \langle F_{ts},F_t\rangle {\bf T}$, where ${\bf N}$ and ${\bf
  T}$ are the unit normal and tangent vectors \cite{Kong08} (then
  $u$ stands for the support function of $F$, the flow supports
  convexity, and finite-time shrinking to a point occurs).

Concerning other nonlinearities, the elliptic M-A equation
 \beq
 \label{GG.13}
 \mbox{$
 ({\rm det} \, D^2 u)^{\frac 1{N+2}} = - \frac 1u, \,\, u<0, \,\, D^2 u > 0
 \,\,\, \mbox{in} \,\,\, \O, \quad u=0 \,\,\, \mbox{on} \,\,\,
 \partial \O,
 $}
  \eeq
where $\O$ is a bounded convex domain in $\ren$, was derived by
Loewner--Nirenberg
 \cite{LN}
in the study of the metric of the form $-\frac 1u \, D^2 u$ ($u$
is then treated as a section of a certain line bundle), and was
proved to admit a unique $C^\infty$ solution; see \cite{CY86} and
earlier references therein. The general Hessian equation has the
form
 \beq
 \label{s**A}
 S_k(D^2 u) = (-u)^p, \,\, u<0 \quad \mbox{in} \,\,\, \O, \quad u=0 \,\,\,
 \mbox{on} \,\,\, \partial \O,
  \eeq
  where $\O$ is a ball in $\ren$, $N \ge 3$,
   and $S_k$ is given by the
   elementary symmetric function
 $$
 \mbox{$
 S_k(D^2u)= \sum_{(1\le i_1<...<i_k \le n)} \l_{i_1}...\l_{i_k},
 $}
 $$
 with  $\{\l_i\}$ being the eigenvalues of the Hessian $D^2 u$ (so
 $k=1$ and $k=n$ correspond to the Laplace and the M-A operators,
 respectively). (\ref{s**A})
  is known to exhibit the {\em critical exponents}
 $$
 \mbox{$
\g(k)=\frac{(N+2)k}{(N-2k)_+},
 $}
 $$
 such that no
 smooth solution $u<0$ exists for $p \ge \g(k)$, and a negative
 radial solution exists for $p \in (0,\g(k))$; see \cite{Tso90}
 (nonexistence is proved by a Pohozaev-type inequality) and \cite{Chou96} for
 extensions. For
 $k=1$, $\g(1)= \frac{N+2}{N-2}$ is the critical Sobolev exponent.
   The nonexistence result for the elliptic equation
 $$
\Delta  u + u^p=0
 $$
is associated with  Pohozaev's classic  inequality \cite{Poh65}.
 The exponents $\g(k)$ above are to be
compared with the critical ones
 $$
\mbox{$ \g(k) = \frac{N+2k}{(N-2k)_+} \quad \mbox{for elliptic
PDEs}$} \quad
-(-\Delta )^k u + |u|^{p-1}u=0,
 $$
 where the
existence-nonexistence results are proved by higher-order
Pohozaev's inequalities \cite{Poh70} applied to general
quasilinear $2k$th-order PDEs. For the M-A equation in a convex
$\O \subset \ren$,
 $$
 {\rm det}\, D^2 u = p(x) g(-u) \quad \mbox{in} \quad \O,
 \quad u=0 \quad \mbox{on} \quad \partial \O,
  $$
  with positive nonincreasing $g$ and $p \in C^\iy$ also positive,
  a criterion of existence of convex negative solution is
  established in \cite{Moh08, Moh09}. Removable singularity theory for elliptic equations (these
include M-A, Hessian, and Weingarten ones)
 $$
 F(D^2 u)= \psi(x,u,D u),
 $$
where e.g., $F(r) \ge F(s) + F(r-s)$ for $r-s>0$ (say, $F$ is
continuous and concave in a set of $N \times N$ real symmetric
matrices) and $F(r) \ge \g ({\rm det}\, r)^{1/N}$, is well
developed; see \cite{Schulz07} for further references and results.
Note that a first removability theorem for the equation
(\ref{DDD1}) in $\re^2$ is due to J\"orgens (1955) \cite{Jor55}.
See the beginning of Section \ref{SectEx} for an extra discussion
of regularity issues for the M-A equations.

As another  standard nowadays direction of elliptic theory, let us
mention recent research on existence of solutions of real
(complex) M-A equations with infinite values (blow-up) on the
boundary $\partial \O$ of a strictly convex bounded smooth domain
$\O \subset \ren$:
 \beq
  \label{in1}
  {\rm det} \, D^2 u(x)=g(x) f(u(x)) \quad \mbox{in} \quad \O,
  \quad u(x) \to \iy \quad \mbox{as}
  \quad x \to \partial \O;
   \eeq
 see \cite{Ivar06, Moh07, Plis08} for history, references, and results.
 Note that, for the semilinear elliptic equations
  \beq
  \label{in2}
  \D u=f(u),
  \eeq
first results were obtained for $f(u)={\mathrm e}^u$ by Bieberbach
in 1916 in 2D \cite{Bieb16} and by Rademacher in 1943 in 3D
followed by  the research for general nonlinearities $f(u)$ by
J.B.~Keller and Osserman in 1957;
   see more historical details in \cite{Moh07}, and
V\'eron \cite{VerMon, VerMon2} and Labutin \cite{Lab03} for a more
complete overview of other results.

\smallskip

 In particular and for instance, the above presentation suggests generalized second-order
M-A parabolic flows (these equations are formulated for $u(x,t)$
being convex or ``almost convex"),
 \beq
 \label{NNN.1Q}
 u_t =  ({\rm det} D^2 u)^{m} \pm (-u)^p,
 \eeq
 with some exponents $m>0$ and $ p \in \re$, that generate many
 {open problems} concerning local existence of
 convex solutions, free-boundary (degeneracy set) propagation,
 extinction, and blow-up singularity patterns, {etc.}
 Other interesting models occur by choosing  the elliptic operator as in
(\ref{GG.11}).
 In radial setting, where (\ref{NNN.1Q}) reduces to
 a 1D quasilinear parabolic PDE,  the interface
 equations and their regularity, moduli of continuity of proper solutions,
  waiting time phenomena, {etc.}, are
 characterized by Sturmian intersection comparison techniques,
 \cite[Ch.~7]{GalGeom}.~For $N>1$, the majority of the problems are
 {open},
  and particular exact  solutions might  be key.
 As we have shown, such models can be considered as natural counterparts of the PME with
 reaction/absorption, and  of thin film (or Cahn--Hilliard-type, $n=0$)
 models, such as (\ref{RD.991}).
On basic properties of hyperbolic M-A equations, see
\cite[Ch.~16]{Tay}. Various parabolic and hyperbolic M-A equations
admitting exact solutions on
 linear invariant  subspaces are described in \cite[Ch.~6]{GSVR}.

Concerning blow-up or extinction behaviour in M-A flows that are
main subject of the present paper and which often are not
well-understood, as a simple illustration, a few types of
singularity formation phenomena occurring on linear subspaces
admitted by such Monge--Amp\`ere operators
 are shown in \cite[\S~6.5]{GSVR}.
 As for (\ref{M1}), these are associated with invariant subspaces
 of the principal M-A operator
 \beq
 \label{He.11}
F_2[u]= {\rm det} D^2u \equiv u_{xx}u_{yy} - (u_{xy})^2 \quad
\mbox{in} \quad \re^2.
 \eeq
 Notice that quadratic polynomials
$p(x)$ occur in the celebrated result of the theory of elliptic
M-A PDEs, establishing that any convex solution of the elliptic
M-A equation
  \beq
  \label{DDD1}
 {\rm det} D^2u=1 \quad \mbox{in} \,\,\, \ren
  \eeq
 is $u(x)=p(x)$. This result  is due to J\"orgens (1954)  for
 $N=2$, Calabi (1958)  for $N =3$, $4$, and $ 5$, and to Pogorelov (1978) for
 any $N \ge 2$
 (see also  \cite{CY86} for a more general result).
 The same conclusion holds for the {\em Hessian
  quotient equation}
 $$
S_{k}(D^2 u)=1 \quad \mbox{for \, $u(x) \le A(1+|x|^2)$ \,
strictly convex},
 $$
   with any $1 \le k < N$ \cite{Bao03}.
  Similarly, if $u(x,t)$ is a smooth solution of
 the parabolic PDE
  $$
  -u_t \, {\rm det} D^2u=1 \quad \mbox{in} \,\,\, \ren \times
  \re_+ \, ,
  $$
  where $u$ is convex in $x$, nonincreasing with $t$, and $u_t$ is
  bounded away from $0$ and $-\infty$, then
 $
u(x,t) = Ct + p(x);
 $
  see \cite{Gut98} for the results and a survey.

\subsection*{On higher-order M-A flows and blow-up}

Fourth and higher-order M-A PDEs have been less well studied,
though some of the equations  correspond to classical geometric
problems, and
several general results have been established. We refer to 
 \cite{Yos00}, where existence was established for the following
 class of
  $m$th-order fully nonlinear equations of the M-A-type
 in $\re^2 \times \re_+$ ($m
\ge 2$)
 \beq
 \label{MA.m1}
 \mbox{$
 u_t = \sum_{(|\mu|=|\nu|=m)} a_{\mu,\nu} D^\mu_x u D^\nu_x u +(\mbox{lower-order
 terms}).
 $}
 \eeq
 Here
  $\mu$ and $\nu$ are multi-indices, and the matrix $\|a_{\mu,\nu}\|$
  satisfies a positivity-type assumption
 for local existence.
In \cite{Yos00}, a Riemann--Hilbert factorization condition
appeared. In \cite{Trud02}, $W^{1,p}$-regularity estimates for
fourth-order M-A equations were derived and  an analogy of the
J\"orgens--Calabi--Pogorelov result for such PDEs was established.
In
 \cite{Li03},  homogeneous fourth-order PDEs for
affine maximal hypersurfaces were studied. Further references can
be found in these papers.

   In order
  to formulate fourth-order Hessian equations in $\re^2$,
let us write
  down the fourth differential of a $C^4$-function $u=u(x,y)$ as a
  {\em quartic form}
   $$
   {\mathrm d}^4 u= u_{xxxx}  {\mathrm d}x^4 + 4 u_{xxxy}  {\mathrm
   d}x^3  {\mathrm d}y + 6 u_{xxyy}  {\mathrm d}x^2 {\mathrm
   d}y^2 + 4 u_{xyyy} {\mathrm d}x {\mathrm d}y^3 + u_{yyyy} {\mathrm
   d}y^4.
 $$
 This gives the {catalecticant determinant} (\ref{F41}).

Using the operator $F_4[u] = {\rm det} \, D^4 u$, it is easy to
construct some formal exact geometric flows. Clearly, $F_4$
preserves the subspace of fourth-degree polynomials
 $
 W_ 7= {\mathcal L}\{1,x^2,xy,y^2,x^4, x^2y^2,y^4\}
 $
 and $F_4: W_7 \to {\mathcal L}\{1\}$.  Therefore,
 the flow
  $ 
  u_t = F_4[u]
  $ 
  is global on $W_7$.
 The basic invariant  subspace of sixth-degree polynomials is
  $
   W= {\mathcal L}\{x^\a y^\b, \,\,0\le \a+\b \le 6\},
 $ 
 on which blow-up may happen {via} a cubic DS.
 Singular patterns also exist for other fourth-order
M-A-type models that are constructed in accordance to their
second-order counterparts. Let us present two other examples with
different types of evolution singularities.

\smallskip



 Let us discuss simple examples on {\em extinction and blow-up} for
 higher-order M-A flows.
    For instance, the M-A equation (cf. (\ref{Loft1}) with the
    known
    extinction behaviour)
 $$
 \mbox{$
 \mbox{$
 u_t = - \frac u{{\rm det} D^4 u} \quad \mbox{in} \,\,\, \re^2
 \times \re_+
 $} $}
  $$
admits   solutions on the subspace $W_4$,
 $
  u(x,y,t)= C_1(t) + C_2(t) x^4 + C_3(t) x^2y^2 + C_4(t) y^4.
  $
  Then $|D^4 u|= 64(36 C_2C_3C_4-C_3^3)$.
The corresponding fourth-order DS yields solutions in separate
variables
 \beq
 \label{HHH.55}
 \mbox{$
  u(x,y,t)= C_1(t)\bigl(1 + A x^4 + B x^2y^2 + C y^4\bigr), 
  $}
  \eeq
 where $A$, $B$, and $C$ are positive constants satisfying $\g = 64(36 ABC-B^3)>0$ by
  the convexity
  assumption on
 initial data. Here  $C_1$ solves the ODE
 $
 C_1'= -
  \frac 1{\g C_1^2}.
 $
  This gives
 finite-time extinction with the rate
   $$
   \mbox{$
   C_1(t)=  \bigl[\frac 3 \g  (T-t)\bigr]^{\frac 13} \to 0 \quad \mbox{as}
   \,\,\,\, t \to T^-.
   $}
   $$
 {Vice versa}, 
  $$
  \mbox{$
  u_t = u \sqrt{{\rm det} D^4 u} \quad \mbox{in} \,\,\, \re^2 \times
  \re_+
  $}
  $$
admits solutions (\ref{HHH.55}) driven by the ODE
 $
 C_1' = \sqrt {\g} \, C_1^{ 5/2}
 $
 with  blow-up,
 $$
 \mbox{$
 C_1(t) =\bigl[ \frac{3 \sqrt \g}2(T-t) \bigr]^{-\frac 23} \to +\infty \quad \mbox{as}
   \,\,\, t \to T^-.
 $}
 $$
 Similar singularity phenomena are traced for the corresponding  hyperbolic
M-A flows  on these   subspaces. A number of typical conclusions
for the second-order Hessian flows can be extended to this
fourth-order, as well as  higher-order, though the well-posedness
of such parabolic or hyperbolic PDEs in classes of ``convex"
functions is a difficult
 {open problem}.

\smallskip

Finally, notice that the homogeneous equation $|D^4 u|=0$ is a
direct sum of two identical copies of the second-order  M-A
equation (see \cite{Fer02})
 \beq
 \label{V.eq1}
 v_{xx}v_{yy} - (v_{xy})^2=0.
 \eeq
 Similarly, the sixth-order equation $|D^6 u|=0$ with the operator
 \beq
 \label{66.11}
{\rm det} D^6 u={\rm det} \, \left[ \begin{matrix}
 u_{60} \,\, u_{51} \,\, u_{42} \,\, u_{33} \cr
  u_{51} \,\, u_{42} \,\, u_{33} \,\, u_{24} \cr
   u_{42} \,\, u_{33} \,\, u_{24} \,\, u_{15} \cr
    u_{33} \,\, u_{24} \,\, u_{15} \,\, u_{06}
    \end{matrix} \right] \quad \big(u_{ij}= u_{x^i y^j}\big) ,
    \eeq
  decouples into three copies of
(\ref{V.eq1}), \cite{Fer02}.~Possibly, this means that some
problems with such higher-order M-A operators are associated with
the second-order ones. In particular, the inhomogeneous equations
$| D^4 u|=1$, $|D^ 6 u|=1$ might  be handled by reduction to
second-order equations, and a result associated with the
J\"orgens--Calabi--Pogorelov theorem might be expected (though
some basics of such PDEs remain obscure).


\end{appendix}

\end{small}

\begin{small}
\end{small}

\end{document}